\documentclass[preprint, 11pt]{elsarticle}
\usepackage{lineno}
\pdfoutput=1

\usepackage{amsmath}
\usepackage{amssymb}
\usepackage{mathrsfs}
\usepackage{color}
\usepackage{xcolor}
\usepackage{graphicx}
\usepackage[percent]{overpic}
\usepackage{url}
\usepackage{caption}
\usepackage{subcaption}
\usepackage{enumerate}
\usepackage{overpic}
\usepackage{amsthm}
\usepackage{moreverb}
\usepackage{multirow}
\usepackage{multicol}
\usepackage[pdftex,colorlinks,bookmarksopen,bookmarksnumbered,citecolor=red,urlcolor=red]{hyperref}
\usepackage{mwe}
\usepackage{graphbox}
\usepackage{textcomp}
\usepackage{tabularx}

\def\x{{\bm x}}

\def\0{\boldsymbol{0}}

\def\cl {\nonumber \\}
\def\el {\nonumber }

\newcommand{\bm}[1]{\mbox{\boldmath{$#1$}}}
\def\div{\nabla\cdot}

\newcommand{\anna}[2][cyan]{{\textcolor{#1}{#2}}}
\newcommand{\rev}[2][red]{{\textcolor{#1}{#2}}}
\newcommand{\lander}[2][orange]{{\textcolor{#1}{#2}}}

\usepackage{verbatim}
\usepackage{graphicx}
\usepackage{epstopdf}

\begin{document}

\begin{frontmatter}


\title{Linear and nonlinear filtering for a two-layer quasi-geostrophic ocean model}





\author[Houston]{Lander Besabe}
\ead{lybesabe@central.uh.edu}

\author[SISSA]{Michele Girfoglio}
\ead{mgirfogl@sissa.it}

\author[Houston]{Annalisa Quaini\corref{mycorrespondingauthor}}
\ead{aquaini@central.uh.edu}
\cortext[mycorrespondingauthor]{Corresponding author}

\author[SISSA]{Gianluigi Rozza}
\ead{grozza@sissa.it}

\address[Houston]{Department of Mathematics, University of Houston, 3551 Cullen Blvd, Houston TX 77204, USA}

\address[SISSA]{SISSA, International School for Advanced Studies, Mathematics Area, mathLab, via Bonomea, Trieste 265 34136, Italy}

\begin{abstract}
Although the two-layer quasi-geostrophic equations (2QGE) are
a simplified model for the dynamics of a stratified, wind-driven ocean, their numerical simulation
is still plagued by the need for high resolution to capture the full spectrum of
turbulent scales. Since such high resolution would lead to unreasonable computational times, 
it is typical to resort to coarse low-resolution meshes 
combined with the so-called eddy viscosity parameterization to account for 
the diffusion mechanisms that are not captured due to mesh under-resolution.
We propose to enable the use of further coarsened meshes
by adding a (linear or nonlinear) differential low-pass to the 2QGE, without
changing the eddy viscosity coefficient. While the linear filter introduces 
constant (additional) artificial viscosity everywhere in the domain, 
the nonlinear filter relies on an indicator function to determine where and how
much artificial viscosity is needed.
Through several numerical results for a double-gyre wind forcing benchmark, we show that with 
the nonlinear filter we obtain accurate results with very coarse meshes, thereby drastically reducing the computational time
(speed up ranging from 30 to 300). 
\end{abstract}


\end{frontmatter}


\section{Introduction}

The simulation of ocean flows is extremely challenging
for several reasons. The first is of course the scale of the problem: the area of an ocean basin is of the order of millions of Km$^2$. The second reason is the nature of the flow itself, which is typically described in terms of 
two non-dimensional numbers: the Reynolds number $Re$ and the Rossby number $Ro$. The Reynolds number is the ratio of inertial forces to viscous forces, while the Rossby number weighs the inertial forces
over the Coriolis forces. Ocean flows characterized by large $Re$ and small $Ro$ are especially challenging 
because they require very fine computational meshes to resolve all the eddy scales. The third reason is connected to the first two: high resolution meshes over very large domains lead to a prohibitive computational cost with
nowadays computational resources. 

In order to contain the computational cost of ocean flow simulations, two avenues are typically followed:
i) to introduce assumptions that reduce the complexity of the model, and ii) to adopt coarse meshes for the simulations. The first important
simplifying assumption comes from the fact that in an ocean basin the vertical length scale is much smaller than horizontal length scales. Thanks to this, one can average the 3D Navier–Stokes equations over the depth and get rid of the vertical dimension, obtaining
the shallow water equations (a 2D problem). If one further assumes hydrostatic and 
geostrophic balance, and small $Ro$ (i.e., flow strongly affected by Coriolis forces), the shallow water equations can be simplified to obtain the quasi-geostrophic equations (QGE).
The QGE, whose name comes from the fact that for $Ro = 0$ one recovers geostrophic flow, 
are the simplest model for wind-driven ocean dynamics.

Since almost all large-scale flows in the ocean (and atmosphere) are in geostrophic balance to leading order,
the QGE are one of the most
commonly used mathematical models
in this context \cite{Majda_Wang_2006,Vallis_2017}, also thanks to their simplicity. However, the QGE
represent the ocean as one layer with uniform depth, density, and temperature, 
while in reality the ocean is a stratified fluid driven from its upper layer by patterns of momentum and buoyancy fluxes \cite{Marshall1997}.
The two-layer quasi-geostrophic equations (2QGE) are an attempt to capture the complexity of stratification 
by adding a second dynamically active layer to include the first baroclinic modes \cite{Marshall1997,Chassignet1998,Berloff1999,DiBattista2001,BERLOFF_KAMENKOVICH_PEDLOSKY_2009}.
While the 2QGE are a more realistic representation 
of ocean dynamics than the QGE, they are also more complex from the mathematical and numerical point of view. 
In fact, the one layer QGE is a two-way coupled system: the kinematic relationship couples the potential vorticity to the stream function and the momentum balance equation has a convective term that involves the stream function. 
The  2QGE system features an additional level of complexity since the kinematic relationships couple the potential vorticity in each layer to the stream functions
of \emph{both} layers, resulting in a tight
coupling between the two layers. 

Even with a model obtained under
simplifying assumptions, like the QGE or 2QGE, a simulation
would require an impractically fine mesh resolution to resolve
the full spectra of turbulence down to the Kolmogorov scale.
To reduce the computational cost to the point that it is manageable, researchers typically resort to coarse low-resolution
meshes and the use of the so-called eddy viscosity  parameterization. This simply means that, instead of using
the actual value of the kinetic viscosity of water 
(of the order of $10^{-6}$ m$^{-2}$/s), the model employs an artificially large viscosity, called eddy viscosity, to compensate for the diffusion
mechanisms that are not captured due to mesh under-resolution. 
The values of the eddy viscosity range from a few thousand (see, e.g., \cite{Bryan1963,Gates1968})
to a few hundred (see, e.g., \cite{Holland1975,Berloff1999,BERLOFF_KAMENKOVICH_PEDLOSKY_2009,Tanaka2010}), depending on the
mesh resolution that one
can afford with the available computational resources.
Of course, this is a rather crude Large Eddy Simulation 
(LES) technique and already a couple of decades ago it was
shown that different eddy viscosity coefficients can result in different dynamics given by quasi-geostrophic models \cite{Berloff1999}. Understanding 
the relationship between eddy viscosity coefficients and mesh resolution has motivated much work in recent years, with the ultimate goal of designing more effective LES methods. Let us focus on the efforts that involve the 2QGE.

LES techniques directly resolve a portion of the flow scales and require a model to account for the remaining (small) scales that are not directly captured due to mesh under-refinement. Since one of the effects of the small (unresolved) scales is to dissipate the turbulent kinetic energy, traditionally,
these models introduce artificial viscosity that can be shown 
to result from an operation of spatial filtering. In \cite{Holm2005}, the Helmholtz filter is applied to the stream function of both layers to obtain a so-called alpha model ($\alpha$ being the filtering radius). Such
alpha model is used to study the mean effects of turbulence on baroclinic instability. In \cite{San2012}, a tridiagonal filter and a differential filter are applied to 
potential vorticity and stream function in both layers, with a closure
approach based on the approximate deconvolution. The work in \cite{San2012} can be seen as an extension to the 2QGE of the LES methodology proposed in \cite{San2011} for the QGE.
In \cite{San_IJMCE2013}, a multi-scale coarse grid projection method
was introduced to smooth out numerical oscillations in the computed solution and 
increase the solution accuracy with a given coarse mesh. 
In \cite{Maulik2016}, an efficient trapezoidal filter is applied to potential vorticity and stream function in both layers, with the Smagorinsky and Leith models
for closure. The difference is that 
the artificial viscosity is proportional to the local strain rate field in the Smagorinsky model, while it is proportional to the gradient of the vorticity field in the Leith model.
It is found that the Smagorinsky model is robust with respect to changes in the grid size, while the
Leith model introduces lower amounts of eddy viscosity
and is less robust to mesh size variations. 
Finally, in  \cite{Berloff2021} local filtering is used to translate the effect of the unresolved scales
into an error-correcting forcing.

In this paper, we propose to apply a linear and nonlinear 
Helmholtz filter (same filter as in \cite{Holm2005} in the linear case) to the potential vorticities (unlike 
\cite{Holm2005} that applies it to the stream functions). In the nonlinear case, our approach uses an indicator function to determine where and how much artificial viscosity is needed. 
The indicator function that we choose makes the closure model similar to the Leith model in \cite{Maulik2016}. Since we use
a different filter than in \cite{Maulik2016}, we find that the 
Leith model is robust with respect to changes in the mesh size
and more accurate than the method based on the linear filter.
Through several numerical results, we show that with the nonlinear {filter,} we obtain accurate results with very coarse meshes, thereby drastically reducing the computational time. The speed up ranges roughly from
30 to 300.

The rest of the paper is structured as follows. Sec.~\ref{sec:pbd} presents the 2QGE model and the application of the linear and nonlinear filter to it. 
Sec.~\ref{sec:discr} describes the methods for space and time discretization. Sec.~\ref{sec:results} reports
and discusses the numerical results. 
Finally, conclusions are drawn in Sec.~\ref{sec:concl}.

\section{Problem definition}
\label{sec:pbd}

\subsection{Two-layer Quasi-Geostrophic Equations}
\label{sec:twoLayerQGE}

Isopycnals are layers within an ocean that are stratified based on their densities. We assume that the part of the ocean under consideration is composed of two isopycnal layers with uniform depth, density, and temperature. Let 
layer 1 be on on top of layer 2. We denote with
$H_1$ and $H_2$ the depths of the layers, both
assumed to be much smaller than the meridional length of the domain $L$. The 2QGE describe the dynamics of these two layers, under some 
simplifying assumptions. 

Two main assumptions for the model are hydrostatic and geostrophic balance.
The second is the so called $\beta$-plane approximation. 
This means that the Coriolis frequency $f$ is linearized as follows: 
$f = f_0 + \beta y$, where $f_0$ is the local rotation rate at $y = 0$, 
which is the center of the basin, and $\beta$ is the gradient of the Coriolis frequency. 
The value of $f_0$ depends on the rotational speed of the earth and the latitude at $y = 0$. 
The $\beta$-plane approximation is equivalent to approximating the Earth (a sphere) with a tangent plane at $y = 0$. Let $\Omega=[x_0, x_f] \times [-L/2, L/2]$ be our computational domain on this plane.

For simplicity, we will work with the 2QGE in 
non-dimensional form. The reader interested in the derivation of these equations from the dimensional form is referred to, e.g., \cite{San2012, Salmon1978, Medjo2000, Fandry1984, Mu1994}.  
In the non-dimensionalization process, some well-known non-dimensional numbers appear. These are the Rossby number $Ro$, the Reynolds number $Re$, and the Froude number $Fr$:
\begin{equation}\label{eq:ro_re_fr}
    Ro = \frac{U}{\beta L^2}, \quad 
    Re = \frac{UL}{\nu}, \quad
    Fr = \frac{f_0^2 U}{g' \beta H},
\end{equation}
where $U$ is a characteristic speed, $\nu$ is the eddy viscosity, and $H = H_1 + H_2$ is the total ocean depth. 
Finally, in \eqref{eq:ro_re_fr}, $g' = g \Delta \rho/\rho_1$, where $g$ is the gravitational constant, $\Delta \rho$ is 
the density jump between the two layers, and $\rho_1$ is the density of the upper layer. The third main assumption in the derivation of the 2QGE is that inertial forces are negligible with respect to the Coriolis and pressure forces, i.e.,  $Ro$ is 
``small''. 

Note that instead of using the actual kinetic viscosity of water (around $10^{-6}$ m$^{-2}$/s), the definition of $Re$ in \eqref{eq:ro_re_fr} uses the eddy viscosity coefficient, which 
can take values of the order of 100-1000 m$^{-2}$/s. This choice is dictated by 
the fact that the scale of an ocean basin is much larger than the effective 
scale for molecular diffusion and so a simulation with a mesh fine enough 
to capture the full spectrum of turbulent scales  
would have a prohibitive cost. To reduce such cost, one uses a coarser mesh
and increases, by trial and error, the value of $\nu$ in order to obtain stable and realistic results. 

We need to introduce two, lesser known, non-dimensional parameters
before stating the 2QGE. These are the lateral eddy viscosity coefficient $A$ and the Ekman bottom layer friction coefficient
$\sigma$:
\begin{equation}\label{eq:a_sig}
    A=\frac{\nu}{\beta L^3}, \quad \sigma = \frac{\gamma}{\beta L},
\end{equation}
where $\gamma$ is a friction coefficient with the bottom of the ocean. Note that $Re = {Ro}/{A}$.
 
Let $q_i$ be the non-dimensional potential vorticity and $\psi_i$ the non-dimensional stream function
of layer $i$, $i= 1,2$. We set $\bm{\Psi}_i=(0, 0, \psi_i)$.
The non-dimensional 2QGE read: find $q_i$ 
and $\psi_i$, for $i= 1,2$, such that 
\begin{align} \label{eq:qge2_1}
    &\frac{\partial q_1}{\partial t} + \div \left(\left(\nabla \times \bm{\Psi}_1 \right)q_1\right) - A\Delta^2 \psi_1 = F_1,\\
    &\frac{\partial q_2}{\partial t} + \div \left(\left(\nabla \times \bm{\Psi}_2 \right)q_2\right) - A\Delta^2 \psi_2 +\sigma\Delta \psi_2 = F_2, \label{eq:qge2_2} \\
    &q_1 = Ro\Delta\psi_1 + y + \frac{Fr}{\delta}\left(\psi_2 - \psi_1 \right), \label{eq:kinematic1}\\
    &q_2 = Ro\Delta\psi_2 + y + \frac{Fr}{1-\delta}\left(\psi_1 - \psi_2 \right),
    \label{eq:kinematic2}
\end{align}
in $\Omega\times(0,T)$, where $F_i$ is the external forcing acting on layer $i$. Here, $(0,T)$ is a time interval of interest. Eq.~\eqref{eq:kinematic1} and \eqref{eq:kinematic2} state the kinematic relationships
between the potential vorticites and the stream functions, with 
$\delta = {H_1}/(H_1+H_2)$ representing the aspect ratio of the layer depths.

Problem \eqref{eq:qge2_1}-\eqref{eq:kinematic2} is rather complex: it is fourth-order in space, nonlinear, and two-way coupled. The two layers are coupled through the last term at the right-hand-side in eq.~\eqref{eq:kinematic1}-\eqref{eq:kinematic2}. Notice that the strength of the coupling is not symmetric, 
as $Fr/\delta$ weights how the dynamics in the bottom layer affect the dynamics in the top layer, while $Fr/(1 - \delta)$
weights the vice versa. So, for given $Fr$ and fixed total depth, the thinner the top layer (i.e., the smaller $\delta$), the stronger the effect of the bottom layer on the top layer and the weaker the effect of the top layer on the bottom layer. 
This is intuitive. We also remark that the dissipative terms in 
\eqref{eq:qge2_1}-\eqref{eq:qge2_2} are multiplied by $A = Ro/Re$ and $\sigma$. Since the definition of $Re$ uses the eddy viscosity, one controls the amount of dissipation in the system by tuning the value of $\nu$.

In order to avoid solving a fourth-order in space problem, 
we use eq.~\eqref{eq:kinematic1} and \eqref{eq:kinematic2} to rewrite eq.~\eqref{eq:qge2_1}-\eqref{eq:qge2_2} as follows:
\begin{align} 
        &\frac{\partial q_1}{\partial t} + \div \left(\left(\nabla \times \bm{\Psi}_1 \right)q_1\right) + \frac{Fr}{Re~\delta}\Delta(\psi_2-\psi_1)
        - \frac{1}{Re}\Delta q_1 = F_1, \label{eq:qge2_1b}\\
      &\frac{\partial q_2}{\partial t} + \div \left(\left(\nabla \times \bm{\Psi}_2 \right)q_2\right) + \frac{Fr}{Re~(1-\delta)}\Delta(\psi_1-\psi_2) - \frac{1}{Re}\Delta q_2 +\sigma\Delta \psi_2 = F_2, \label{eq:qge2_2b}
\end{align}
in $\Omega\times(0,T)$.

Obviously, we need to supplement the model with boundary and initial conditions. Following \cite{San2012}, we prescribe free-slip and impermeability boundary conditions
\begin{align} 
    \psi_i = 0 &\quad \mbox{on }\partial\Omega\times(t_0,T), \label{qge-bdry1}\\
    q_i = y &\quad \mbox{on }\partial\Omega\times(t_0,T), \label{qge-bdry2}
\end{align}
for $i=1,2$, and start the system from a rest state:
\begin{align} 
    q_i(x,y,t_0) = y &\quad \mbox{in }\Omega. \label{qge-initial}
\end{align}

\subsection{Low-pass differential spatial filter}
\label{subsec:filter}

Once the value of the eddy viscosity $\nu$ is set, 
a Direct Numerical Simulation (DNS) of the 2QGE requires a mesh whose size $h$ is smaller than the Munk scale:
\begin{align}
\delta_M = L \, \sqrt[3]{\dfrac{Ro}{Re}}, \label{eq:munk}
\end{align}
where $L$ is a characteristic length. 
We are slightly abusing of terminology, since, strictly speaking, DNS
refers to a simulation that uses a mesh with 
size $h < \delta_M$ where $\delta_M$ is computed with $Re$ based
on the value of molecular diffusivity. However, a real DNS of ocean dynamics is far beyond the computing resources available today. 
We propose to couple the 2QGE with a differential filter to enable the use of even coarser meshes than required by the Munk scale computed with $Re$ based on the eddy viscosity. 
The purpose of the differential filter is to model the effects of the scales that are not resolved due to mesh under-refinement.

Let $\alpha$ be a filtering radius and $a(\cdot)$ 
a scalar function, called {indicator function}, with the following properties:
\begin{align*}
    0<a(\cdot)\leq1 & \mbox{ for any } (\x, t)\in\Omega\times[t_0,T];\\
    a(\cdot)\simeq 0 & \mbox{ in regions where the flow field needs no regularization;}\\
    a(\cdot)\simeq 1 & \mbox{ in regions where the flow field needs } O(\alpha) \mbox{ regularization.}
\end{align*}
The extension of the filtering strategy proposed in \cite{Girfoglio_JCAM2023}
for the single layer QGE to the 2QGE reads: 
\begin{align}
        &\frac{\partial q_1}{\partial t} + \nabla \cdot \left(\left(\nabla \times \bm{\Psi}_1 \right)q_1\right) + \frac{Fr}{Re~\delta}\Delta(\psi_2-\psi_1)
        - \frac{1}{Re}\Delta q_1 = F_1, \label{eq:coupled-sys-1}\\
        &-\alpha^2 \nabla \cdot \left( a(q_1)\nabla \bar{q}_1\right) + \bar{q}_1=q_1, \label{eq:coupled-sys-2} \\ 
        & Ro\Delta\psi_1 + y + \frac{Fr}{\delta}\left(\psi_2 - \psi_1 \right) = \bar{q}_1 , \label{eq:coupled-sys-3}\\
        &\frac{\partial q_2}{\partial t} + \nabla \cdot \left(\left(\nabla \times \bm{\Psi}_2 \right)q_2\right) + \frac{Fr}{Re~(1-\delta)}\Delta(\psi_1-\psi_2)- \frac{1}{Re}\Delta q_2 +\sigma\Delta \psi_2 = F_2,\label{eq:coupled-sys-4}\\
        &-\alpha^2 \nabla \cdot \left( a(q_2)\nabla \bar{q}_2\right) + \bar{q}_2=q_2, \label{eq:coupled-sys-5} \\ 
        &Ro\Delta\psi_2 + y + \frac{Fr}{1-\delta}\left(\psi_1 - \psi_2 \right) = \bar{q}_2, \label{eq:coupled-sys-6}
\end{align}
in $\Omega\times(0,T)$, where $\bar{q}_1$ and $\bar{q}_2$ are the filtered potential vorticities 
for layer 1 and 2, respectively. Notice that, while 
\eqref{eq:coupled-sys-1} and \eqref{eq:coupled-sys-4}
are the same as \eqref{eq:qge2_1b}-\eqref{eq:qge2_2b}, eq.~\eqref{eq:coupled-sys-3} replaces $q_1$ in \eqref{eq:kinematic1}
with $\bar{q}_1$ and, analogously, 
eq.~\eqref{eq:coupled-sys-6} replaces $q_2$ in \eqref{eq:kinematic2}
with $\bar{q}_2$. 
It is trivial to see that 
as $\alpha \rightarrow 0$, 
$\bar{q}_1 \rightarrow q_1$ and
$\bar{q}_2 \rightarrow q_2$
and system \eqref{eq:coupled-sys-1}-\eqref{eq:coupled-sys-6} tends to 
system \eqref{eq:kinematic1}-\eqref{eq:qge2_2b}. For $\alpha > 0$,
the differential filters
\eqref{eq:coupled-sys-2} and \eqref{eq:coupled-sys-5}
leverage an elliptic
operator that acts as a spatial filter by damping the spurious and nonphysical oscillations exhibited by the numerical solution on coarse meshes. The price that we pay to have a physical solution with a coarse grid is the addition of two equations, i.e., eq.~\eqref{eq:coupled-sys-2} and \eqref{eq:coupled-sys-5}.
As we shall see in Sec.~\ref{sec:results}, this additional cost is marginal with respect
to the gain in computational time allowed by 
the use of a coarse grid. 

Filter problems \eqref{eq:coupled-sys-2} and \eqref{eq:coupled-sys-5} are linear or nonlinear depending on the choice of the 
indicator function. 
If $a(\cdot) = 1$, eq.~\eqref{eq:coupled-sys-2} and \eqref{eq:coupled-sys-5} are linear and one recovers an extension to the 2QGE of the so-called  BV-$\alpha$ model for the one-layer case \cite{Monteiro2015, Nadiga2001, Girfoglio_JCAM2023, Holm2003, Monteiro2014}.
We will call this 2QG-$\alpha$ model. 
While $a(\cdot)=1$ is a convenient choice,
it does not really act as an indicator function since it introduces the same
amount of regularization everywhere in the domain. If  $a(\cdot)$ is actually a function of its input, i.e., not a constant, then it can be selective in introducing
regularization. Following \cite{Girfoglio_JCAM2023}, we will consider
\begin{equation} \label{eq:nl-ind-fxn}
    a(q) = \frac{|\nabla q|}{||\nabla q||_\infty},
\end{equation}
which makes eq.~\eqref{eq:coupled-sys-2} and \eqref{eq:coupled-sys-5} nonlinear. We will refer to the problem \eqref{eq:coupled-sys-1}-\eqref{eq:coupled-sys-6} with 
indicator function \eqref{eq:nl-ind-fxn} 
as the 2QG-NL-$\alpha$ model, which is an extension to the 2QGE of the BV-NL-$\alpha$ for the single layer QGE \cite{Girfoglio_JCAM2023}. 

We supplement problem \eqref{eq:coupled-sys-1}-\eqref{eq:coupled-sys-6} with initial
data \eqref{qge-initial} and boundary conditions \eqref{qge-bdry1}-\eqref{qge-bdry2}, plus the following additional boundary conditions:
\begin{align} 
    \bar{q}_i = y, &\quad \mbox{on }\partial\Omega\times(0,T). \el
\end{align}
for $i=1,2$.

\section{The discretized problem}\label{sec:discr}

We will start with the time discretization of 
problem~\eqref{eq:coupled-sys-1}-\eqref{eq:coupled-sys-6} and then
present the space discretization of the time-discrete problem. 

\subsection{Time discretization}
\label{subsec:time-discrete}

We divide the interval $(0,T)$ into $N_T$ subintervals of width $\Delta t = \frac{T}{N_T}$ and denote $t^n = n\Delta t$, with $n = 0, 1, \dots, N_T$. For any given quantity $f$, its approximation at a specific time point $t^n$ is denoted with $f^n$. 

Problem~\eqref{eq:coupled-sys-1}-\eqref{eq:coupled-sys-6} discretized in time by a Backward Differentiation Formula of order 1 (BDF1) reads: for $i = 1,2$,
given $\left( q_i^0,\bar{q}_i^0, \psi_i^0 \right)$, 
find $(q^{n+1}_i, \bar{q}^{n+1}_i,\psi^{n+1}_i)$, for $n\geq0$, such that
\begin{align}
        &\frac{1}{\Delta t}q_1^{n+1} + \div\left(\left(\nabla \times \bm{\Psi}_1^{n}\right)q_1^{n+1}\right) + \frac{Fr}{Re~\delta} \Delta \left( \psi_2^{n+1} - \psi_1^{n+1} \right) - \frac{1}{Re}\Delta q_1^{n+1} = b_1^{n+1}, \label{gqe2-discretized-time}
     \\
    & -\alpha^2 \nabla \cdot \left( a_1^{n+1}\nabla \bar{q}_1^{n+1}\right) + \bar{q}_1^{n+1} = q_1^{n+1} \\
    & Ro\Delta \psi_1^{n+1} + y + \frac{Fr}{\delta}\left(\psi_2^{n+1} - \psi_1^{n+1}\right) = \bar{q}_1^{n+1}, \\
    & \frac{1}{\Delta t}q_2^{n+1} + \div\left(\left(\nabla \times \bm{\Psi}_2^{n}\right)q_2^{n+1}\right) + \frac{Fr}{Re~(1-\delta)}\Delta \left( \psi_1^{n+1} - \psi_2^{n+1} \right) \cl
    & \quad - \frac{1}{Re}\Delta q_2^{n+1} + \sigma \Delta \psi_2^{n+1}= b_2^{n+1} \label{gqe2-discretized-time1} \\
    & -\alpha^2 \nabla \cdot \left( a_2^{n+1}\nabla \bar{q}_2^{n+1}\right) + \bar{q}_2^{n+1} = q_2^{n+1} \\
    & Ro\Delta \psi_2^{n+1} + y + \frac{Fr}{1-\delta}\left(\psi_1^{n+1} - \psi_2^{n+1}\right) = \bar{q}_2^{n+1} 
    \label{gqe2-discretized-time2}
\end{align}
where $b_i^{n+1} = F_i + q_i^n/\Delta t$ and $a_i^{n+1} = a(q_i^{n+1})$. Note that we have linearized the convective terms in \eqref{gqe2-discretized-time} and \eqref{gqe2-discretized-time1} with a first-order extrapolation. However, 
if one uses a nonlinear filter, problem \eqref{gqe2-discretized-time}-\eqref{gqe2-discretized-time2} remains nonlinear. 
Additionally, it is coupled. Hence, one needs to be careful in devising a
solution algorithm that contains the computational cost
while achieving a desired level of accuracy. 


To maximize computational efficiency, we propose a 
segregated algorithm that at the same time linearizes the problem when using a nonlinear filter. At time step $t^{n+1}$, 
given $(q_i^n, \bar{q}_i^n, \psi_i^n)$, for $i=1,2$, proceed with the following steps:

\begin{itemize}
    \item[-] Step 1: find the potential vorticity of the top layer $q_1^{n+1}$ such that
    \begin{align}\label{eq:seg-alg1}
        \frac{1}{\Delta t}q_1^{n+1} + \div\left(\left(\nabla \times \bm{\Psi}_1^{n}\right)q_1^{n+1}\right) - \frac{1}{Re}\Delta q_1^{n+1} = b_1^{n+1} 
        - \frac{Fr}{Re~\delta} \Delta \left( \psi_2^{n} - \psi_1^{n} \right).
    \end{align}
    \item[-] Step 2: find the filtered potential vorticity for the top layer $\bar{q}_1^{n+1}$  such that
    \begin{equation} \label{eq:seg-alg2}
        -\alpha^2 \nabla \cdot \left( a_1^{n+1}\nabla \bar{q}_1^{n+1}\right) + \bar{q}_1^{n+1} = q_1^{n+1}.
    \end{equation}
    \item[-] Step 3: find the stream function of the top layer $\psi_1^{n+1}$ such that:
    \begin{equation} \label{sf-1-discreet}
         Ro\Delta \psi_1^{n+1} + y - \frac{Fr}{\delta} \psi_1^{n+1} = \bar{q}_1^{n+1} - \frac{Fr}{\delta}\psi_2^{n}.
    \end{equation}
    \item[-] Step 4: find the potential vorticity  of the bottom layer $q_2^{n+1}$ such that:
    \begin{align}
        &\frac{1}{\Delta t}q_2^{n+1} + \div\left(\left(\nabla \times \bm{\Psi}_2^{n}\right)q_2^{n+1}\right) - \frac{1}{Re}\Delta q_2^{n+1}= b_2^{n+1} -\sigma \Delta \psi_2^{n} \cl
        & \quad - \frac{Fr}{Re~(1 - \delta)} \Delta \left( \psi_1^{n+1} - \psi_2^{n} \right).
    \end{align}   
    \item[-] Step 5: find the filtered potential vorticity for the bottom layer $\bar{q}_2^{n+1}$ such that:
    \begin{equation} \label{eq:seg-alg5}
        -\alpha^2 \nabla \cdot \left( a_2^{n+1}\nabla \bar{q}_2^{n+1}\right) + \bar{q}_2^{n+1} = q_2^{n+1}.
    \end{equation}
    \item[-] Step 6: find the stream function of the bottom layer $\psi_2^{n+1}$ such that:
    \begin{equation}\label{eq:seg-alg6}
         Ro\Delta \psi_2^{n+1} + y - \frac{Fr}{1-\delta} \psi_2^{n+1} = \bar{q}_2^{n+1} - \frac{Fr}{1-\delta}\psi_1^{n+1} .
    \end{equation}
\end{itemize}

\subsection{Space discretization}
\label{subsec:space-discrete}

For the space discretization, we adopt a Finite Volume approximation that is derived directly from the integral form of the governing equations.
Computational domain $\Omega$ is divided into a total number of $N_C$ control volumes $\Omega_k$, $k=1,\dots,N_C$. Let  \textbf{A}$_j$ be the surface vector of each face of the control volume, 
with $j = 1, \dots, M$.

The integral form of \eqref{eq:seg-alg1} in each control volume $\Omega_k$ is:
\begin{align} \label{eq:int-form1}
    \frac{1}{\Delta t}\int_{\Omega_k}q_1^{n+1}\,d\Omega + \int_{\Omega_k} \div\left(\left(\nabla \times \bm{\Psi}_1^{n}\right) q_1^{n+1}\right) \,d\Omega - \frac{1}{Re}\int_{\Omega_k}\Delta q_1^{n+1}\,d\Omega \cl 
    = \int_{\Omega_k}b_1^{n+1}\,d\Omega - \frac{Fr}{Re~\delta}\int_{\Omega_k}\Delta \left( \psi_2^n - \psi_1^n \right)\, d\Omega_k.
\end{align}
By applying the Gauss-divergence theorem to \eqref{eq:int-form1}, we obtain
\begin{align} \label{gauss-div}
    \frac{1}{\Delta t}\int_{\Omega_k}q_1^{n+1}\,d\Omega + \int_{\partial \Omega_k} \left(\left(\nabla \times \bm{\Psi}_1^{n}\right) q_1^{n+1}\right) \cdot d\textbf{A} - \frac{1}{Re}\int_{\Omega_k}\nabla q_1^{n+1} \cdot d\textbf{A} \cl 
    = \int_{\Omega_k}b_1^{n+1}\,d\Omega - \frac{Fr}{Re~\delta}\int_{\partial\Omega_k}\nabla \left( \psi_2^n - \psi_1^n \right) \cdot d\textbf{A}.
\end{align}
The discretized form of eq.~\eqref{gauss-div}, divided by the control volume 
$\Omega_k$, can be written as:
\begin{align} \label{eq:disc-q1}
        \frac{1}{\Delta t}q_{1,k}^{n+1} + \sum_j\varphi_{1,j}^{n}q_{1,k}^{n+1,j} - \frac{1}{Re}\sum_j \left(\nabla q_{1,k}^{n+1}\right)_j\cdot \textbf{A}_j = b_{1,k}^{n+1}
        - \frac{Fr}{Re~\delta}\sum_j \left(\nabla(\psi_{2,k}^n - \psi_{1,k}^n)\right)_j\cdot \textbf{A}_j
\end{align}
with 
\begin{equation}\label{eq:flux1}
\varphi_{1,j}^n=\left(\nabla \times \bm{\Psi}_{1,j}^{n}\right)\cdot\textbf{A}_j,  
\end{equation}
where $q_{1,k}^{n+1}$ and $\psi_{1,k}^{n}$ represent the average top layer potential vorticity and average stream function in the control volume $\Omega_k$, $q_{1,k}^{n+1,j}$ is the top layer potential vorticity associated to the centroid of the $j$-th face and normalized by the control volume $\Omega_k$, and
$b_{1,k}^{n+1}$ denotes the average discrete forcing on the top layer in $\Omega_k$. 
The convective term \eqref{eq:flux1}, as well as its counterpart for the bottom layer \eqref{eq:flux2}, 
are computed by linear interpolation from neighboring cells using central difference, which is a second-order method. The Laplacian terms for both layers are approximated with a second order central differencing scheme. 
More details on the treatment of the convective and diffusive terms can be found in \cite{Girfoglio2019, GirfoglioPSIZETA}.

By following the same procedure for 
eqs. \eqref{eq:seg-alg2}-\eqref{eq:seg-alg6} and using similat notation, we get:
\begin{align}
    & -\alpha^2 \sum_j a_{1,k}^{n+1} \left(\nabla \bar{q}_{1,k}^{n+1}\right)_j\cdot \textbf{A}_j + \bar{q}_{1,k}^{n+1}=q_{1,k}^{n+1},\\
    & Ro\sum_j\left(\nabla \psi_{1,k}^{n+1}\right)_j\cdot\textbf{A}_j + y_k + \frac{Fr}{\delta}\left(\psi_{2,k}^n-\psi_{1,k}^{n+1}\right) = \bar{q}_{1,k}^{n+1},\\
    & \frac{1}{\Delta t}q_{2,k}^{n+1} + \sum_j\varphi_{2,j}^{n}q_{2,k}^{n+1,j} - \frac{1}{Re}\sum_j \left(\nabla q_{2,k}^{n+1}\right)_j\cdot \textbf{A}_j = b_{2,k}^{n+1} \cl 
    & \quad + \left(\frac{Fr}{Re~(1-\delta)} - \sigma\right) \sum_j (\nabla\psi_{2,k}^n)_j\cdot \textbf{A}_j - \frac{Fr}{Re~(1-\delta)}\sum_j \left(\nabla\psi_{1,k}^{n+1}\right)_j\cdot \textbf{A}_j,\\
    & -\alpha^2 \sum_j a_{1,k}^{n+1} \left(\nabla \bar{q}_{1,k}^{n+1}\right)_j\cdot \textbf{A}_j + \bar{q}_{1,k}^{n+1}=q_{1,k}^{n+1},\\
     & Ro\sum_j\left(\nabla \psi_{2,k}^{n+1}\right)_j\cdot\textbf{A}_j + y_k - \frac{Fr}{1-\delta}\psi_{2,k}^{n+1}= \bar{q}_{2,k}^{n+1} - \frac{Fr}{1-\delta}\psi_{1,k}^{n+1}, \label{eq:disc-psi2}
\end{align}
where $y_k$ is the vertical coordinate of the centroid and
\begin{equation}\label{eq:flux2}
\varphi_{2,j}^n=\left(\nabla \times \bm{\Psi}_{2,j}^{n}\right)\cdot\textbf{A}_j.  
\end{equation}

We implement scheme \eqref{eq:disc-q1}-\eqref{eq:flux2} in the open-source software package Geophysical and Environmental Applications (GEA) \cite{GEA, GirfoglioFVCA10}, which 
builds upon C++ finite volume library OpenFOAM\textsuperscript{\textregistered} \cite{Weller1998}. See \cite{Girfoglio_AIP2023, Clinco_2023, Hajisharifi2024} for more details on GEA.

\section{Numerical results} \label{sec:results}

This section presents several numerical experiments to validate our solver for the 2QGE and show the improvements in terms of accuracy and computational efficiency made possible by the
filtering approach presented in Sec.~\ref{sec:discr}.
The validation of the 2QGE solver using a manufactured solution
is discussed in Sec.~\ref{sec:validation}. 
Our investigation on the performance of the 2QG-$\alpha$ and 2QG-NL-$\alpha$ models with a
double-gyre wind forcing benchmark is presented in 
Sec.~\ref{sec:numtests}. 

\subsection{Validation of our solver for the 2QGE} \label{sec:validation}

To validate our segregated Finite Volume solver for the 2QGE, we consider a manufactured steady-state solution solution with varying parameters in domain $\Omega = [-0.5,0.5]^2$. 
Hence, the meridional length $L$ is 1.
We set the stream functions to:
\begin{align} \label{eq:ss-psi1}
    \psi_1 = A_1(x^2 - 0.25)(y^2 - 0.25), \\
    \psi_1 = A_2(x^2 - 0.25)(y^2 - 0.25), \label{eq:ss-psi2}
\end{align}
where $A_1, A_2 > 0$. By plugging these into \eqref{eq:kinematic1}-\eqref{eq:kinematic2}, we obtain the potential vorticities:
\begin{align} \label{eq:ss-q1}
    & q_1 = 2A_1Ro(x^2 + y^2 - 0.5) + y + \frac{Fr}{\delta}(A_2-A_1)(x^2 - 0.25)(y^2 - 0.25), \\
    & q_2 = 2A_2Ro(x^2 + y^2 - 0.5) + y + \frac{Fr}{1-\delta}(A_1-A_2)(x^2 - 0.25)(y^2 - 0.25), \label{eq:ss-q2}
\end{align}
with which we can find 
the forcing terms $F_1$ and $F_2$ 
through \eqref{eq:qge2_1b}-\eqref{eq:qge2_2b}: 
\begin{align}
    F_1 &= 8RoA_1^2xy(x^2 - 0.25) - (8RoA_1^2xy + 2A_1x)(y^2 - 0.25) - \frac{8A_1Ro}{Re},\\
    F_2 &= 8RoA_2^2xy(x^2 - 0.25) - (8RoA_2^2xy + 2A_2x - 2A_2\sigma)(y^2 - 0.25) - \frac{8A_2Ro}{Re}.
\end{align}
To obtain this manufactured solution, we impose the following boundary conditions
\begin{align}
    &\psi_i = 0 && \text{on } \partial \Omega, \\
    &q_i = 2A_iRo(x^2-0.25)-0.5 && \text {on } y = -0.5, \\
    &q_i = 2A_iRo(x^2-0.25)+0.5 && \text {on } y = 0.5, \\
    &q_i = 2A_iRo(y^2-0.25)+y && \text {on } x = -0.5, 0.5,
\end{align}
for $i=1,2$.

Despite the fact that we are using a steady-state solution for simplicity, 
we use a time-dependent solver for the 2QGE, i.e., the solver described 
in Sec \ref{subsec:space-discrete} in the limit of $\alpha = 0$.
As initial condition, we take $q_i(0,x,y)$ given by \eqref{eq:ss-q1} and \eqref{eq:ss-q2}, i.e., 
we give the solver the exact solution and checks that it maintains it.


We set $A_1 = 1$, $A_2 = 2$, $Fr = 0.1$, $\delta = 0.2$, $\sigma = 0$, and vary the remaining two critical parameters: the Reynolds number
and the Rossby number. We proceed as follows: first we set $Ro = 1$ and consider $Re = 10, 100, 1000$, 
then we  set $Re = 1$ and consider $Ro = 0.1, 0.01,0.001$. 
Note that the case $Ro = 1$-$Re = 10$
has the same Munk scale \eqref{eq:munk} as the 
$Re = 1$-$Ro = 0.1$, i.e., $\delta_M \approx 0.46$.
Similarly, the cases $Ro = 1$-$Re = 100$
and $Re = 1$-$Ro = 0.01$ have the same 
Munk scale $\delta_M \approx 0.21$, while the cases
$Ro = 1$-$Re = 1000$ and $Re = 1$-$Ro = 0.001$ have 
Munk scale $\delta_M \approx 0.1$. We consider a set
of meshes that are all fine enough to resolve
these Munk scales. Specifically, we consider mesh sizes $h = 1/32,1/64, 1/128, 1/256$. 
Hence, there is no need for the filter. 
We further remark that, while the test cases
from the first and second set are match by Munk scale, the first set of tests have varying Kolmogorov scale, since $Re$ varies. Defined as
\begin{equation}\label{eq:eta-re}
\eta=Re^{-3/4}L,
\end{equation}
the Kolmogorov scale of the second set of tests
is simply 1. For the first set of tests, it takes the values of 
$\eta = 0.178$ for $Re = 10$, 
$\eta = 0.032$ for $Re = 100$, and
$\eta = 0.006$ for $Re = 1000$.
All the meshes under consideration satisfy the requirement $h < \eta$
for $Re = 10, 100$, while only the mesh with $h = 1/256$ satisfies it for $Re = 1000$. Given that we are working with a manufactured solution, not with a realistic application, we do not expect this to be a problem. 

Tables \ref{table:ms-poly-Re10}-\ref{table:ms-poly-Re1000} report the relative $L^2$ errors for each variable for fixed $Ro$ and varying $Re$ and the corresponding convergence rates. 
Note that 
the rate of convergence is either 2.00 or very close to it
for all the variables. Since we use second order approximations for the spatial derivatives, this is the expected rate and we observe that it does not
deteriorate even when the Reynolds number is increased.

\begin{table}[htb!]
    \centering
    \begin{tabular}{|c|c|c|c|c|c|c|c|c|}
       \hline
       {mesh}  & \multicolumn{2}{c|}{$\psi_{1}$} & \multicolumn{2}{c|}{$\psi_{2}$} & \multicolumn{2}{c|}{$q_{1}$} & \multicolumn{2}{c|}{$q_{2}$} \\ \cline{2-9} 
       size  & error & rate & error & rate & error & rate & error & rate \\ \hline
        $1/32$ & 1.99E-03 & & 1.99E-03 & & 6.17E-04 & & 6.90E-04 & \\  \hline
        $1/64$ & 4.97E-04 & 2.00 & 4.97E-04 & 2.00 & 1.54E-04 & 2.00 & 1.72E-04 & 2.00 \\ \hline
        $1/128$ & 1.24E-04 & 2.00 & 1.24E-04 & 2.00 & 3.86E-05 & 2.00 & 4.31E-05 & 2.00 \\ \hline
        $1/256$ & 3.12E-05 & 1.99 & 3.08E-05 & 2.01 & 9.74E-06 & 1.99 & 1.06E-05 & 2.03 \\ \hline
    \end{tabular}
    \caption{Validation: relative $L^2$ error for each variable and corresponding rates of convergence for different mesh resolutions for $Ro=1$ and $Re=10$.}
    \label{table:ms-poly-Re10}
\end{table}

\begin{table}[htb!]
    \centering
    \begin{tabular}{|c|c|c|c|c|c|c|c|c|}
       \hline
       {mesh}  & \multicolumn{2}{c|}{$\psi_{1}$} & \multicolumn{2}{c|}{$\psi_{2}$} & \multicolumn{2}{c|}{$q_{1}$} & \multicolumn{2}{c|}{$q_{2}$} \\ \cline{2-9} 
        size & error & rate & error & rate & error & rate & error & rate \\ \hline
        $1/32$ & 2.06E-03 & & 2.03E-03 & & 7.36E-04 & & 7.56E-04 & \\  \hline
        $1/64$ & 5.15E-04 & 2.00 & 5.07E-04 & 2.00 & 1.84E-04 & 2.00 & 1.89E-04 & 2.00 \\ \hline
        $1/128$ & 1.29E-04 & 2.00 & 1.27E-04 & 2.00 & 4.59E-05 & 2.00 & 4.72E-05 & 2.00 \\ \hline
        $1/256$ & 3.20E-05 & 2.01 & 3.14E-05 & 2.01 & 1.14E-05 & 2.01 & 1.16E-05 & 2.02 \\ \hline
    \end{tabular}
    \caption{Validation: relative $L^2$ error for each variable and corresponding rates of convergence for different mesh resolutions for $Ro=1$ and $Re=100$.}
    \label{table:ms-poly-Re100}
\end{table}

\begin{table}[htb!]
    \centering
    \begin{tabular}{|c|c|c|c|c|c|c|c|c|}
       \hline
       {mesh}  & \multicolumn{2}{c|}{$\psi_{1}$} & \multicolumn{2}{c|}{$\psi_{2}$} & \multicolumn{2}{c|}{$q_{1}$} & \multicolumn{2}{c|}{$q_{2}$} \\ \cline{2-9} 
        size & error & rate & error & rate & error & rate & error & rate \\ \hline
        $1/32$ & 2.28E-03 & & 2.09E-03 & & 1.02E-03 & & 8.50E-04 & \\  \hline
        $1/64$ & 5.70E-04 & 2.00 & 5.22E-04 & 2.00 & 2.55E-04 & 2.01 & 2.12E-04 & 2.00 \\ \hline
        $1/128$ & 1.42E-04 & 2.00 & 1.31E-04 & 2.00 & 6.35E-05 & 2.01 & 5.30E-05 & 2.00 \\ \hline
        $1/256$ & 3.51E-05 & 2.02 & 3.25E-05 & 2.01 & 1.57E-05 & 2.02 & 1.32E-05 & 2.01 \\ \hline
    \end{tabular}
    \caption{Validation: relative $L^2$ error for each variable and corresponding rates of convergence for different mesh resolutions for $Ro=1$ and $Re=1000$.}
    \label{table:ms-poly-Re1000}
\end{table}

\begin{figure}[htb!]
\centering
\begin{tabular}{ccccc}
       & $Re = 10$ & $Re = 100$ & $Re = 1000$ & \\
$\psi_1$  & \includegraphics[align=c,scale = 0.3]{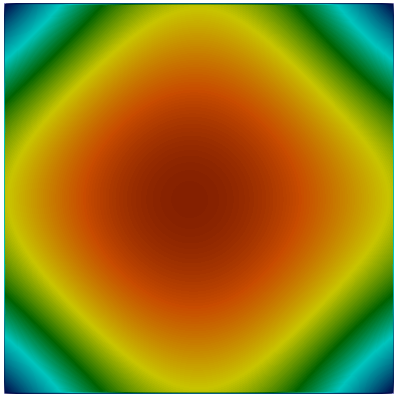} & \includegraphics[align=c,scale = 0.3]{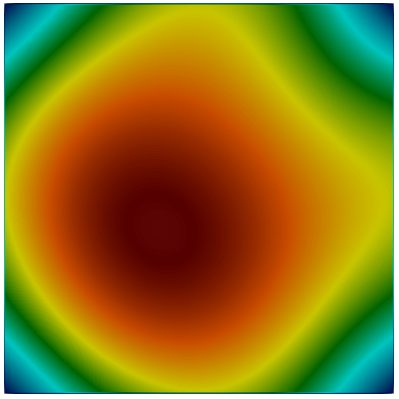} & \includegraphics[align=c,scale = 0.3]{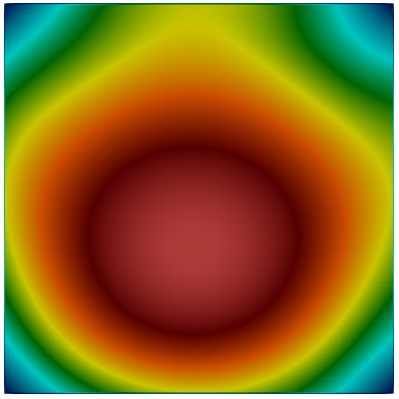} & \hspace{-0.4cm}\includegraphics[align=c,scale = 0.3]{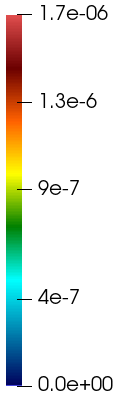} \\
$\psi_2$  & \includegraphics[align=c,scale = 0.3]{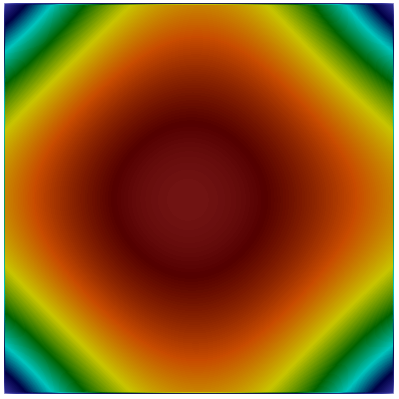} & \includegraphics[align=c,scale = 0.3]{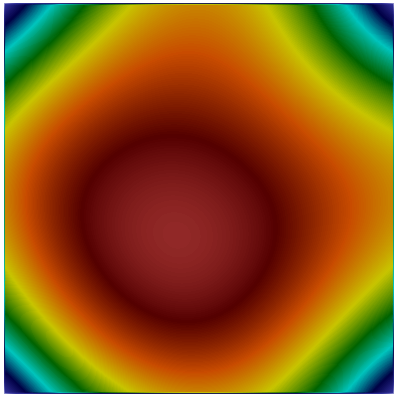} & \includegraphics[align=c,scale = 0.3]{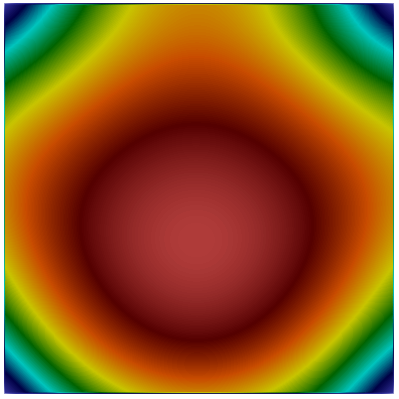} & \hspace{-0.4cm}\includegraphics[align=c,scale = 0.3]{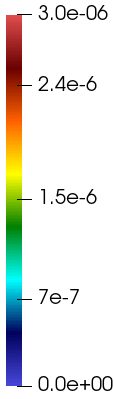}
\end{tabular}
\caption{Validation: absolute error of the stream functions $\psi_1$ (top row) and $\psi_2$ (bottom row) for $Re = 10$ (first column), $Re = 100$ (second column), and $Re = 1000$ (third column). In each case, $Ro = 1$. The results were obtained with mesh size $1/256$.}
\label{fig:poly_re_error1}
\end{figure}

Fig.~\ref{fig:poly_re_error1} illustrates
the difference in absolute value between exact and computed stream functions for $Ro=1$ and varying 
$Re$ for mesh size $h = 1/256$. 
We observe that the maximum error is of order 1E-06.

\begin{figure}[htb!]
\centering
\begin{tabular}{ccccc}
       & $Ro = 0.1$ & $Ro = 0.01$ & $Ro = 0.001$ & \\
$\psi_1$ & \includegraphics[align=c,scale = 0.3]{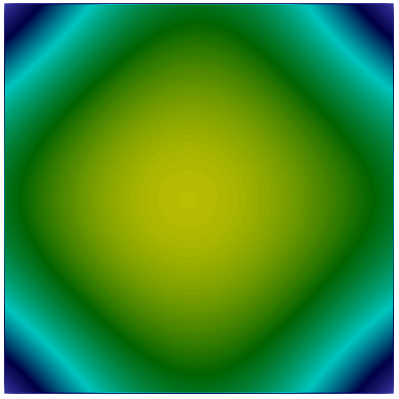} & \includegraphics[align=c,scale = 0.3]{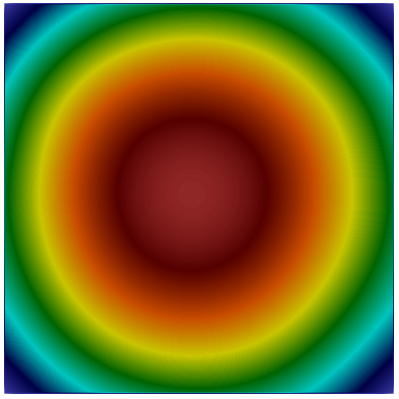} & \includegraphics[align=c,scale = 0.3]{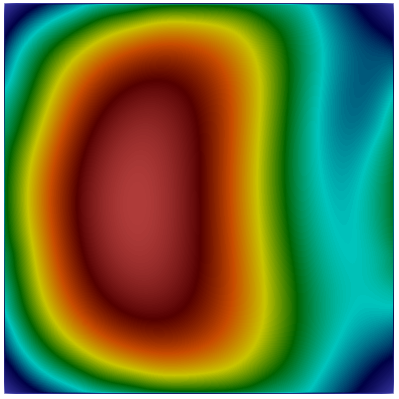} & \hspace{-0.4cm}\includegraphics[align=c,scale = 0.3]{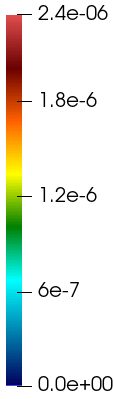} \\
$\psi_2$ & \includegraphics[align=c,scale = 0.3]{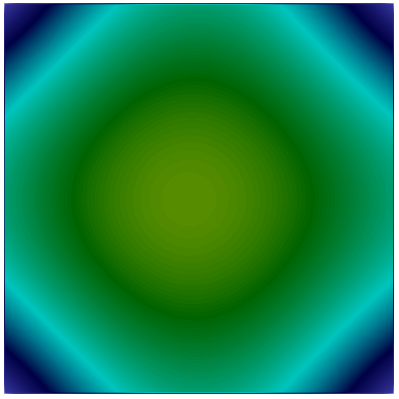} & \includegraphics[align=c,scale = 0.3]{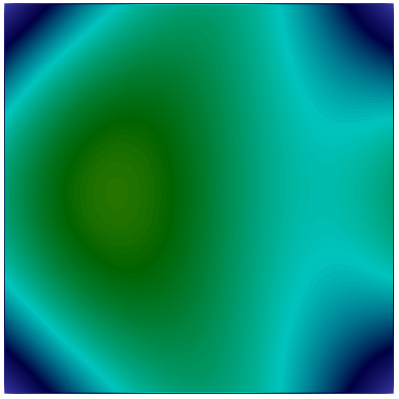} & \includegraphics[align=c,scale = 0.3]{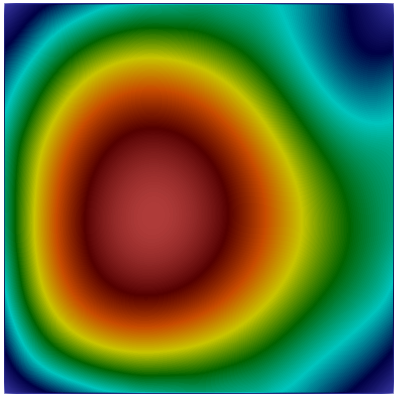} & \hspace{-0.4cm}\includegraphics[align=c,scale = 0.3]{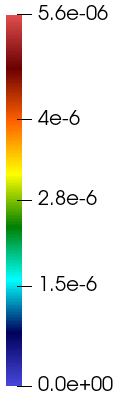}
\end{tabular}
\caption{Validation: absolute error of the stream functions $\psi_1$ (top row) and $\psi_2$ (bottom row) for $Ro = 0.1$ (first column), $Ro = 0.01$ (second column), and $Ro = 0.001$ (third column). In each case, $Re = 1$. The results were obtained with mesh size $h=1/256$.}
\label{fig:poly_re_error}
\end{figure}

Tables \ref{table:ms-poly-Ro0.1}--\ref{table:ms-poly-Ro0.001} report the relative $L^2$ errors for each variable for fixed $Re$ and varying $Ro$ and the corresponding convergence rates. 
We see that
the rate of convergence is 2.00 or very close to it
for all the variables when $Re=1$ and $Ro=0.1$ and for most
variables in the other two cases. We observe
a slight deterioration of the convergence rate for 
$\psi_1$ when $Ro=0.01$ 
and for $q_1$ when $Ro=0.001$. Although this set of cases has matching Munk scales with the previous set,
this is somewhat expected. In fact, it is known that 
the lower $Ro$ cases are more challenging than the higher $Re$
cases, the Munk scale being equal. The reason is the following: a relatively small Rossby number yields a sharp Western boundary
layer, i.e., a convection-dominated regime that makes the problem more challenging. 
See, e.g., \cite{QGE-review} and references therein.
The effect of this sharp Western boundary layer is evident in Fig.~\ref{fig:poly_re_error}, which shows the absolute errors of the stream functions $\psi_1$ and $\psi_2$ 
obtained with mesh $h = 1/256$ for $Re = 1$ and varying $Ro$.
We see that the red region, corresponding to larger errors, moves
towards West as $Ro$ decreases. 
Additionally, we note that 
the maximum absolute difference is slightly higher than in the first set of tests, indicating that these cases are more challenging. Nonetheless, the errors in  Fig.~\ref{fig:poly_re_error1} and \ref{fig:poly_re_error}
have the same order of magnitude, i.e., 1E-06.

\begin{table}[htb!]
    \centering
    \begin{tabular}{|c|c|c|c|c|c|c|c|c|}
       \hline
       {mesh}  & \multicolumn{2}{c|}{$\psi_{1}$} & \multicolumn{2}{c|}{$\psi_{2}$} & \multicolumn{2}{c|}{$q_{1}$} & \multicolumn{2}{|c|}{$q_{2}$} \\ \cline{2-9} 
        size & error & rate & error & rate & error & rate & error & rate \\ \hline
        $1/32$ & 1.99E-03 & & 1.99E-03 & & 6.39E-05 & & 3.29E-04 & \\  \hline
        $1/64$ & 4.97E-04 & 2.00 & 4.97E-04 & 2.00 & 1.60E-05 & 2.00 & 8.22E-04 & 2.00 \\ \hline
        $1/128$ & 1.24E-04 & 2.00 & 1.24E-04 & 2.00 & 3.86E-06 & 2.05 & 2.06E-05 & 2.00 \\ \hline
        $1/256$ & 3.06E-05 & 2.02 & 3.13E-05 & 1.99 & 9.09E-07 & 2.09 & 5.28E-06 & 1.97 \\ \hline
    \end{tabular}
    \caption{Validation: relative $L^2$ error for each variable and corresponding rates of convergence for different mesh resolutions for $Re=1$ and $Ro=0.1$.}
    \label{table:ms-poly-Ro0.1}
\end{table}

\begin{table}[htb!]
    \centering
    \begin{tabular}{|c|c|c|c|c|c|c|c|c|}
       \hline
       {mesh}  & \multicolumn{2}{c|}{$\psi_{1}$} & \multicolumn{2}{c|}{$\psi_{2}$} & \multicolumn{2}{c|}{$q_{1}$} & \multicolumn{2}{c|}{$q_{2}$} \\ \cline{2-9} 
        size & error & rate & error & rate & error & rate & error & rate \\ \hline
        $1/32$ & 2.09E-03 & & 2.10E-03 & & 1.02E-04 & & 7.09E-05 & \\  \hline
        $1/64$ & 5.24E-04 & 2.00 & 5.22E-04 & 2.01 & 2.53E-05 & 2.02 & 1.75E-04 & 2.01 \\ \hline
        $1/128$ & 1.35E-04 & 1.95 & 1.29E-04 & 2.01 & 5.83E-06 & 2.12 & 4.26E-06 & 2.04 \\ \hline
        $1/256$ & 4.36E-05 & 1.63 & 2.79E-05 & 2.21 & 6.40E-07 & 3.19 & 7.42E-07 & 2.52 \\ \hline
    \end{tabular}
    \caption{Validation: relative $L^2$ error for each variable and corresponding rates of convergence for different mesh resolutions for $Re=1$ and $Ro=0.01$.}
    \label{table:ms-poly-Ro0.01}
\end{table}

\begin{table}[htb!]
    \centering
    \begin{tabular}{|c|c|c|c|c|c|c|c|c|}
       \hline
       {mesh}  & \multicolumn{2}{c|}{$\psi_{1}$} & \multicolumn{2}{c|}{$\psi_{2}$} & \multicolumn{2}{c|}{$q_{1}$} & \multicolumn{2}{c|}{$q_{2}$} \\ \cline{2-9} 
        size & error & rate & error & rate & error & rate & error & rate \\ \hline
        $1/32$ & 3.28E-03 & & 3.27E-03 & & 1.83E-04 & & 5.63E-05 & \\  \hline
        $1/64$ & 8.12E-04 & 2.01 & 8.16E-04 & 2.00 & 4.62E-05 & 1.99 & 1.41E-04 & 1.99 \\ \hline
        $1/128$ & 2.20E-04 & 1.88 & 2.03E-04 & 2.01 & 1.05E-05 & 2.14 & 3.28E-06 & 2.11 \\ \hline
        $1/256$ & 4.30E-05 & 2.36 & 4.85E-05 & 2.07 & 3.11E-06 & 1.75 & 9.42E-07 & 1.80 \\ \hline
    \end{tabular}
    \caption{Validation: relative $L^2$ error for each variable and corresponding rates of convergence for different mesh resolutions for $Re=1$ and $Ro=0.001$.}
    \label{table:ms-poly-Ro0.001}
\end{table}

\subsection{Assessment of the linear and nonlinear filtering} \label{sec:numtests}

We consider the extension of the so-called double-wind gyre forcing test, a classical benchmark widely used to assess solvers for geophysical flow problems  \cite{Nadiga2001, Holm2003, Greathbatch2000, Monteiro2014, Monteiro2015, Girfoglio_JCAM2023, QGE-review, San2015, Girfoglio2023}.

We take the computational domain $\Omega$ is $[0,1]\times[-1,1]$, thus $L = 2$, 
and set the forcing terms $F_1 = \sin(\pi y)$ and $F_2 = 0$. Taking inspiration from \cite{Girfoglio_JCAM2023, San2012, Girfoglio2023}, we consider two sets of parameters:
\begin{itemize}
    \item[-] Case 1: $Ro = 0.001$, $Re = 450$, $Fr = 0.1$, $\sigma = 0.005$, and $\delta = 0.5$.
    \item[-] Case 2: $Ro = 0.001$, $Re = 450$, $Fr = 0.1$, $\sigma = 0.01$, and $\delta = 0.1$.
\end{itemize}
These two cases have the same $A = Ro/Re = 0.22$E-5 and hence the same amount of ``physical'' diffusivity and same Munk scale $\delta_M =0.026$. %
The two cases differ in the friction with the bottom of the ocean (less friction in case 1) and the relative depths of the layers (thinner top layer in case 2).

For the DNS of the 2QGE, we use a mesh with size $h = 1/256$ (about \lander{3} times smaller than the Munk scale), denoted as $256\times512$. To test the effect of the filtering, we will use coarse meshes $64\times128$, $32\times64$, $16\times32$ and $8\times16$. The time interval of interest is $[0, 100]$, with time step $\Delta t =$ 2.5E-05 \cite{Girfoglio_JCAM2023}. 
When the filter is used, 
we set {$\alpha = O(h)$} \cite{Holm2003, Girfoglio_JCAM2023, Girfoglio2023}, where $h$ is the mesh size, as is typically done when filtering stabilization is used.
See, e.g., \cite{BQV} and reference therein for more details on this.

\begin{figure}[htb!]
\centering
    \begin{subfigure}{0.2\textwidth}
         \centering
         \includegraphics[width=\textwidth]{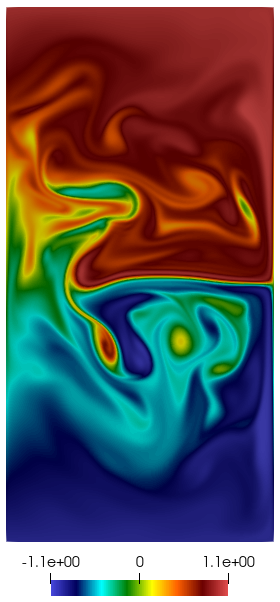}
         \caption{\scriptsize{$q_1$ at $t=20$}}
     \end{subfigure}
    \begin{subfigure}{0.2\textwidth}
         \centering
         \includegraphics[width=\textwidth]{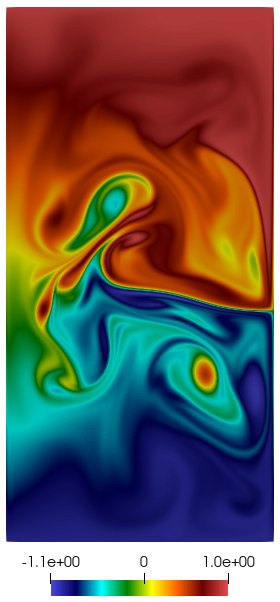}
         \caption{\scriptsize{$q_1$ at $t=21$}}
     \end{subfigure}
    \begin{subfigure}{0.2\textwidth}
         \centering
         \includegraphics[width=\textwidth]{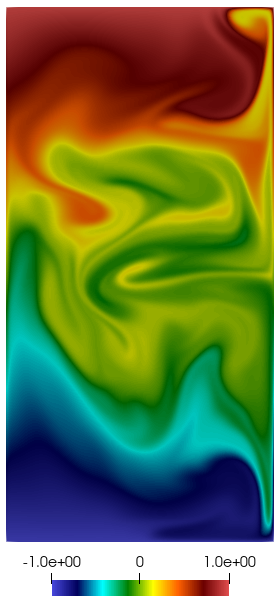}
         \caption{\scriptsize{$q_2$ at $t=20$}}
     \end{subfigure}
    \begin{subfigure}{0.2\textwidth}
         \centering
         \includegraphics[width=\textwidth]{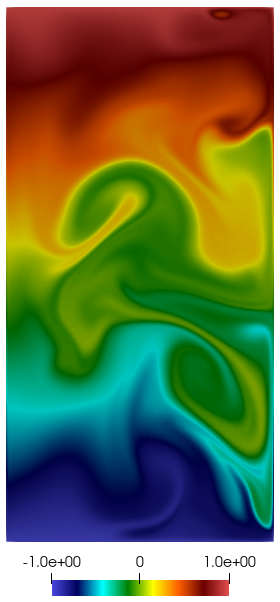}
         \caption{\scriptsize{$q_2$ at $t=21$}}
     \end{subfigure}
\caption{Case 1:  potential vorticity of the top  and bottom layers at time $t=20,21$ computed by the DNS. 
}
\label{fig:PV-instantaneous}
\end{figure}

The quantities of interest for this benchmark are the time-averaged stream functions $\widetilde{\psi}_i$ and potential vorticities $\widetilde{q}_i$, $i = 1, 2$, over the time period $[20,100]$. The reason for this choice is that these quantities 
reach statistically a steady state, while
the instantaneous variables display a highly convective behavior that makes the comparison between DNS solutions and LES solutions harder. Fig.~\ref{fig:PV-instantaneous} 
illustrates this almost chaotic behavior with
an example of $q_1$ and $q_2$ computed by the DNS for case 1 at
time $ t= 20, 21$. Despite being only 1 time unit apart, the vorticity
fields in both layers show significant differences.
Additionally, we will 
track the evolution of  the enstrophy $\mathcal{E}_i$ of each layer:
\begin{equation}\label{eq:entrophy}
    \mathcal{E}_i (t) = \int_{\Omega}q_i^2\,d\Omega, \quad i  = 1,2.
\end{equation}

\subsubsection{Numerical results for case 1} \label{sec:experiment1}

\begin{figure}
\centering
    \begin{overpic}[percent,width=0.16\textwidth,grid=false]{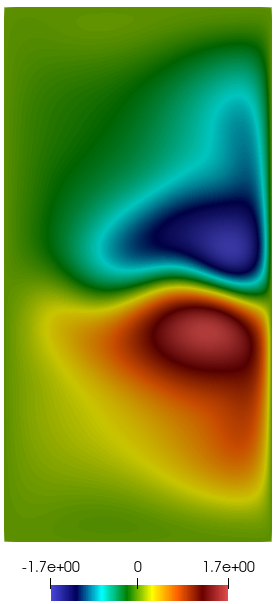}
    \put(5,92){\textcolor{white}{\scriptsize{DNS, $256 \times 512$}}}
    \end{overpic}
    \begin{overpic}[percent,width=0.162\textwidth,grid=false]{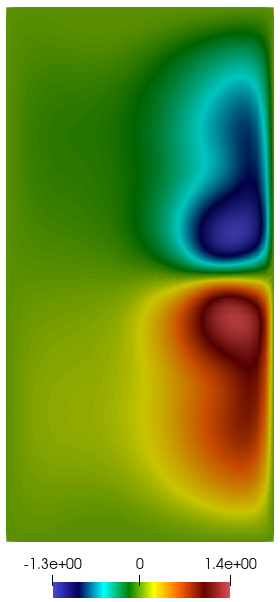}
    \put(6,92){\textcolor{white}{\scriptsize{2QGE, $32 \times 64$}}}
    \end{overpic}
    \begin{overpic}[percent,width=0.162\textwidth,grid=false]{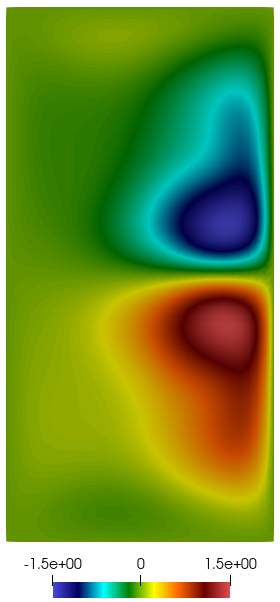}
    \put(5,92){\textcolor{white}{\scriptsize{2QG-$\alpha$, $32 \times 64$}}}
    \end{overpic}
        \begin{overpic}[percent,width=0.162\textwidth,grid=false]{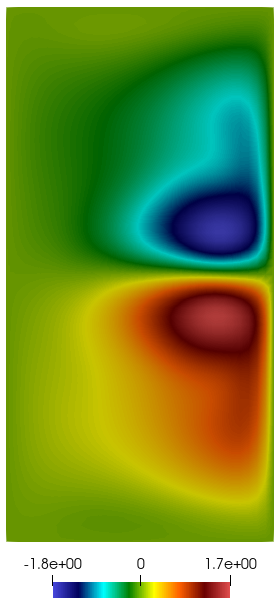}
    \put(1,92){\textcolor{white}{\scriptsize{2QG-NL-$\alpha$, $32 \times 64$}}}
    \end{overpic}\\
    \begin{overpic}[percent,width=0.16\textwidth,grid=false]{img/mf-d0.5/DNS-psi1mean0.5.png}
    \put(5,92){\textcolor{white}{\scriptsize{DNS, $256 \times 512$}}}
    \end{overpic}
    \begin{overpic}[percent,width=0.162\textwidth,grid=false]{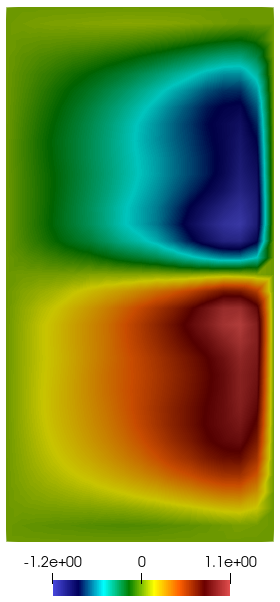}
    \put(6,92){\textcolor{white}{\scriptsize{2QGE, $16 \times 32$}}}
    \end{overpic}
        \begin{overpic}[percent,width=0.162\textwidth,grid=false]{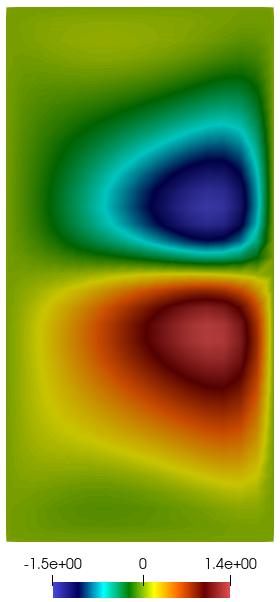}
    \put(5,92){\textcolor{white}{\scriptsize{2QG-$\alpha$, $16 \times 32$}}}
    \end{overpic}
        \begin{overpic}[percent,width=0.162\textwidth,grid=false]{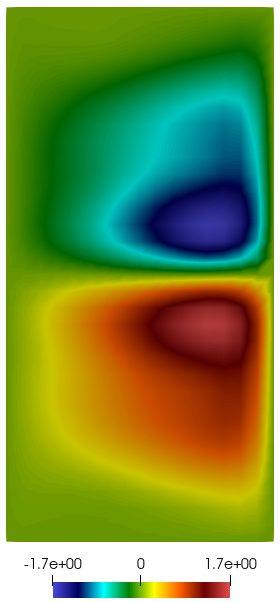}
    \put(1,92){\textcolor{white}{\scriptsize{2QG-NL-$\alpha$, $16 \times 32$}}}
       \end{overpic}\\
   \begin{overpic}[percent,width=0.16\textwidth,grid=false]{img/mf-d0.5/DNS-psi1mean0.5.png}
    \put(5,92){\textcolor{white}{\scriptsize{DNS, $256 \times 512$}}}
    \end{overpic}
        \begin{overpic}[percent,width=0.162\textwidth,grid=false]{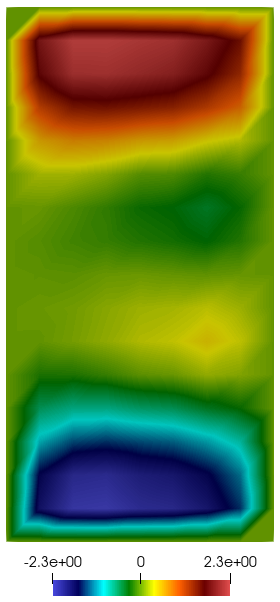}
    \put(6,92){\textcolor{white}{\scriptsize{2QGE, $8 \times 16$}}}
    \end{overpic}
         \begin{overpic}[percent,width=0.162\textwidth,grid=false]{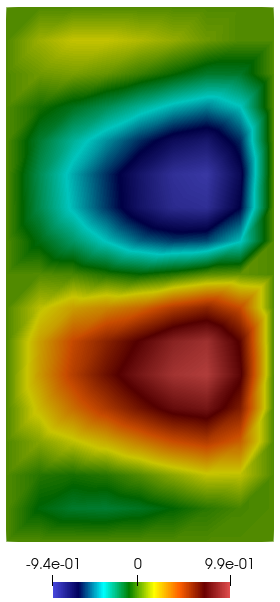}
    \put(5,92){\textcolor{white}{\scriptsize{2QG-$\alpha$, $8 \times 16$}}}
    \end{overpic}
    \begin{overpic}[percent,width=0.162\textwidth,grid=false]{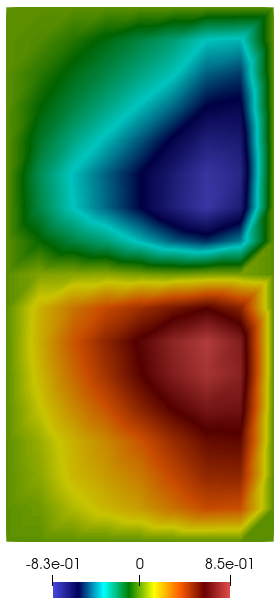}
    \put(2,92){\textcolor{white}{\scriptsize{2QG-NL-$\alpha$, $8 \times 16$}}}
       \end{overpic}
\caption{{Case 1:} Time-averaged stream function of the top layer $\widetilde{\psi}_1$ computed by DNS (first column) and 2QGE with no filtering (second column), 2QG-$\alpha$ (third column), and 2QG-NL-$\alpha$ (fourth column) with different coarse meshes. Note that the color bar may differ
from one panel to the other.}
\label{fig:psi1-case1-filtering}
\end{figure}

\begin{figure}
\centering
    \begin{overpic}[percent,width=0.16\textwidth,grid=false]{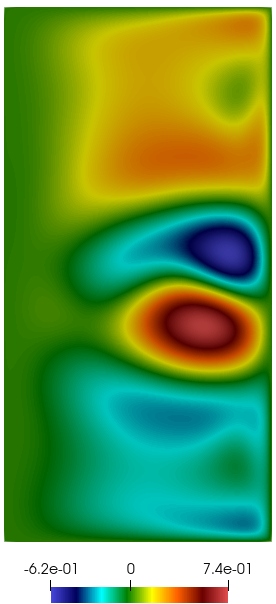}
    \put(5,92){\textcolor{white}{\scriptsize{DNS, $256 \times 512$}}}
    \end{overpic}
    \begin{overpic}[percent,width=0.162\textwidth,grid=false]{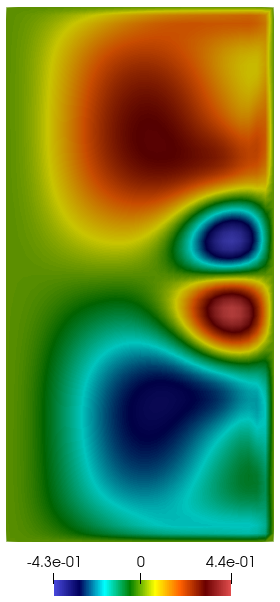}
    \put(6,92){\textcolor{white}{\scriptsize{2QGE, $32 \times 64$}}}
    \end{overpic}
    \begin{overpic}[percent,width=0.162\textwidth,grid=false]{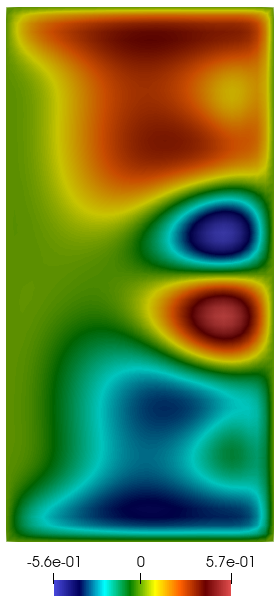}
    \put(5,92){\textcolor{white}{\scriptsize{2QG-$\alpha$, $32 \times 64$}}}
    \end{overpic}
        \begin{overpic}[percent,width=0.162\textwidth,grid=false]{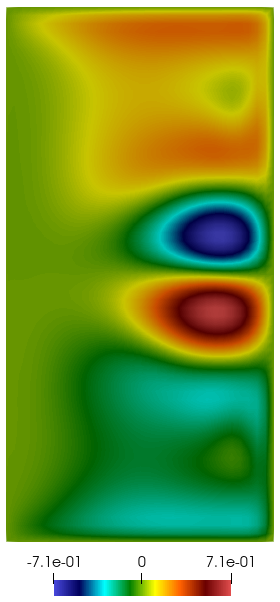}
    \put(1,92){\textcolor{white}{\scriptsize{2QG-NL-$\alpha$, $32 \times 64$}}}
    \end{overpic}\\
    \begin{overpic}[percent,width=0.16\textwidth,grid=false]{img/mf-d0.5/DNS-psi2mean0.5.png}
    \put(5,92){\textcolor{white}{\scriptsize{DNS, $256 \times 512$}}}
    \end{overpic}
    \begin{overpic}[percent,width=0.162\textwidth,grid=false]{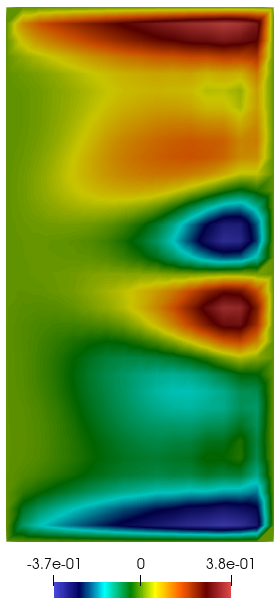}
    \put(6,92){\textcolor{white}{\scriptsize{2QGE, $16 \times 32$}}}
    \end{overpic}
        \begin{overpic}[percent,width=0.162\textwidth,grid=false]{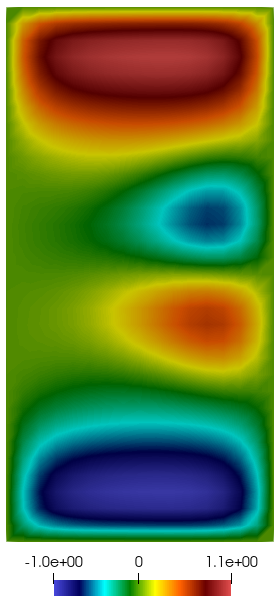}
    \put(5,92){\textcolor{white}{\scriptsize{2QG-$\alpha$, $16 \times 32$}}}
    \end{overpic}
        \begin{overpic}[percent,width=0.162\textwidth,grid=false]{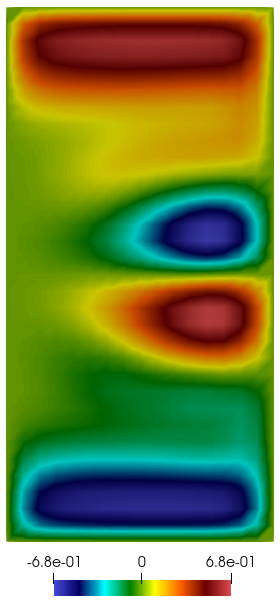}
    \put(1,92){\textcolor{white}{\scriptsize{2QG-NL-$\alpha$, $16 \times 32$}}}
       \end{overpic}\\
   \begin{overpic}[percent,width=0.16\textwidth,grid=false]{img/mf-d0.5/DNS-psi2mean0.5.png}
    \put(5,92){\textcolor{white}{\scriptsize{DNS, $256 \times 512$}}}
    \end{overpic}
        \begin{overpic}[percent,width=0.162\textwidth,grid=false]{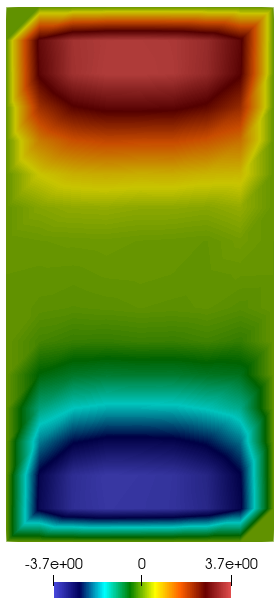}
    \put(6,92){\textcolor{white}{\scriptsize{2QGE, $8 \times 16$}}}
    \end{overpic}
         \begin{overpic}[percent,width=0.162\textwidth,grid=false]{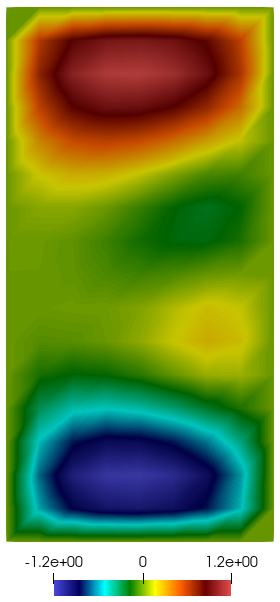}
    \put(5,92){\textcolor{white}{\scriptsize{2QG-$\alpha$, $8 \times 16$}}}
    \end{overpic}
    \begin{overpic}[percent,width=0.162\textwidth,grid=false]{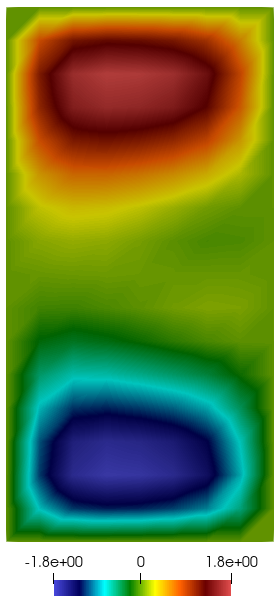}
    \put(2,92){\textcolor{white}{\scriptsize{2QG-NL-$\alpha$, $8 \times 16$}}}
       \end{overpic}
\caption{{Case 1:} Time-averaged stream function of the bottom layer $\widetilde{\psi}_2$ computed by DNS (first column) and 2QGE with no filtering (second column), 2QG-$\alpha$ (third column), and 2QG-NL-$\alpha$ (fourth column) with different coarse meshes. Note that the color bar may differ
from one panel to the other.}
\label{fig:psi2-case1-filtering}
\end{figure}

\begin{figure}
\centering
\begin{overpic}[percent,width=0.16\textwidth,grid=false]{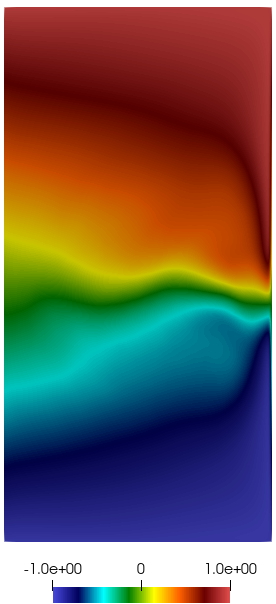}
    \put(5,92){\textcolor{white}{\scriptsize{DNS, $256 \times 512$}}}
    \end{overpic}
    \begin{overpic}[percent,width=0.162\textwidth,grid=false]{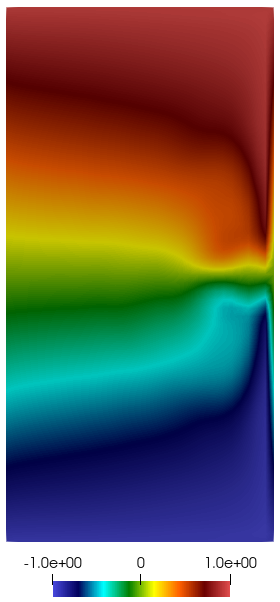}
    \put(6,92){\textcolor{white}{\scriptsize{2QGE, $32 \times 64$}}}
    \end{overpic}
    \begin{overpic}[percent,width=0.162\textwidth,grid=false]{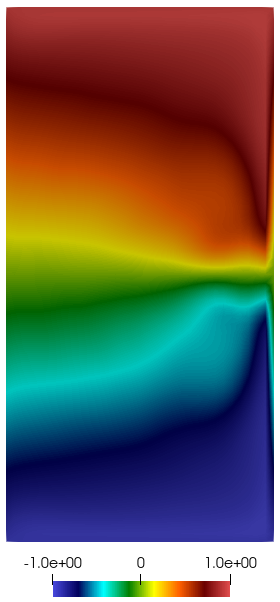}
    \put(5,92){\textcolor{white}{\scriptsize{2QG-$\alpha$, $32 \times 64$}}}
    \end{overpic}
        \begin{overpic}[percent,width=0.162\textwidth,grid=false]{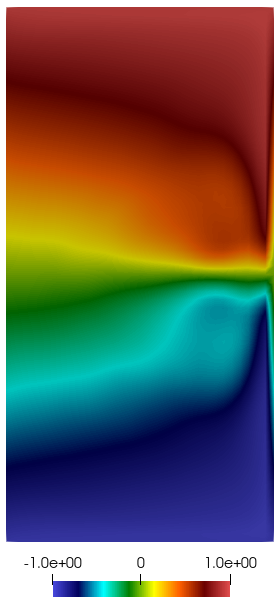}
    \put(1,92){\textcolor{white}{\scriptsize{2QG-NL-$\alpha$, $32 \times 64$}}}
    \end{overpic}\\
\begin{overpic}[percent,width=0.16\textwidth,grid=false]{img/mf-d0.5/DNS-q1mean0.5.png}
    \put(5,92){\textcolor{white}{\scriptsize{DNS, $256 \times 512$}}}
    \end{overpic}
    \begin{overpic}[percent,width=0.162\textwidth,grid=false]{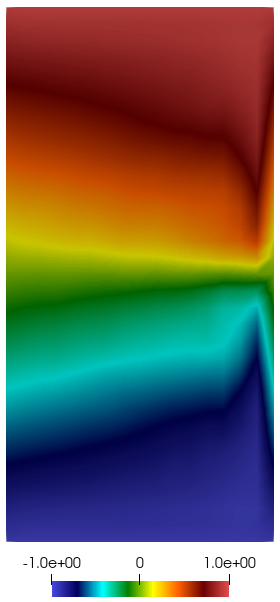}
    \put(6,92){\textcolor{white}{\scriptsize{2QGE, $16 \times 32$}}}
    \end{overpic}
        \begin{overpic}[percent,width=0.162\textwidth,grid=false]{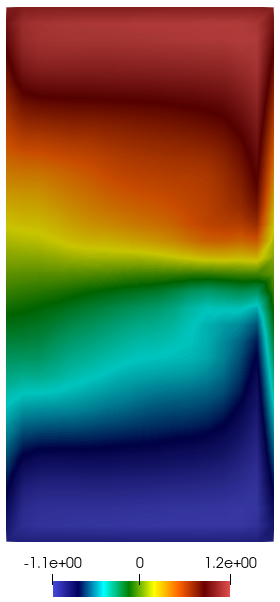}
    \put(5,92){\textcolor{white}{\scriptsize{2QG-$\alpha$, $16 \times 32$}}}
    \end{overpic}
        \begin{overpic}[percent,width=0.162\textwidth,grid=false]{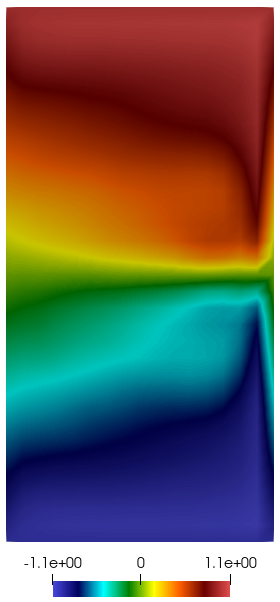}
    \put(1,92){\textcolor{white}{\scriptsize{2QG-NL-$\alpha$, $16 \times 32$}}}
       \end{overpic}\\
   \begin{overpic}[percent,width=0.16\textwidth,grid=false]{img/mf-d0.5/DNS-q1mean0.5.png}
    \put(5,92){\textcolor{white}{\scriptsize{DNS, $256 \times 512$}}}
    \end{overpic}
        \begin{overpic}[percent,width=0.162\textwidth,grid=false]{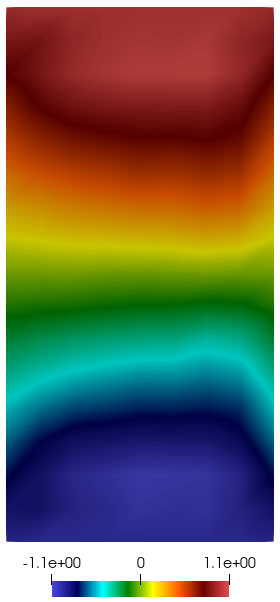}
    \put(6,92){\textcolor{white}{\scriptsize{2QGE, $8 \times 16$}}}
    \end{overpic}
         \begin{overpic}[percent,width=0.162\textwidth,grid=false]{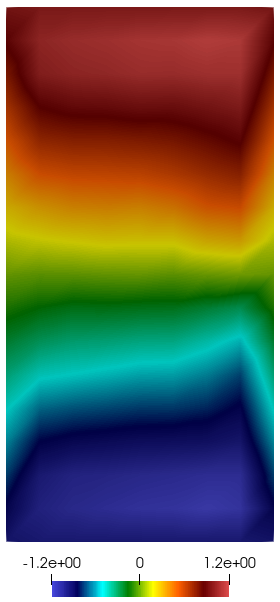}
    \put(5,92){\textcolor{white}{\scriptsize{2QG-$\alpha$, $8 \times 16$}}}
    \end{overpic}
    \begin{overpic}[percent,width=0.162\textwidth,grid=false]{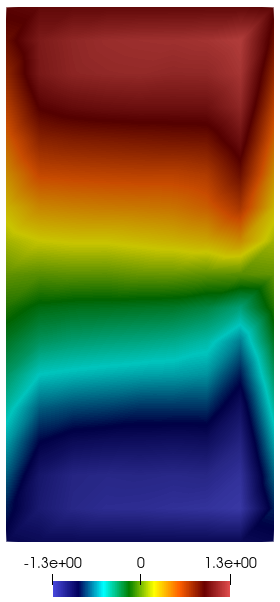}
    \put(2,92){\textcolor{white}{\scriptsize{2QG-NL-$\alpha$, $8 \times 16$}}}
       \end{overpic}
\caption{{Case 1:} Time-averaged potential vorticity of the top layer $\widetilde{q}_1$ computed by DNS (first column) and 2QGE with no filtering (second column), 2QG-$\alpha$ (third column), and 2QG-NL-$\alpha$ (fourth column) with different coarse meshes. Note that the color bar may differ
from one panel to the other.}
\label{fig:q1-case1-filtering}
\end{figure}

\begin{figure}
\centering
    \begin{overpic}[percent,width=0.16\textwidth,grid=false]{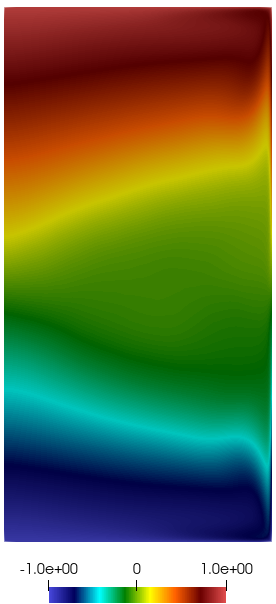}
    \put(5,92){\textcolor{white}{\scriptsize{DNS, $256 \times 512$}}}
    \end{overpic}
    \begin{overpic}[percent,width=0.162\textwidth,grid=false]{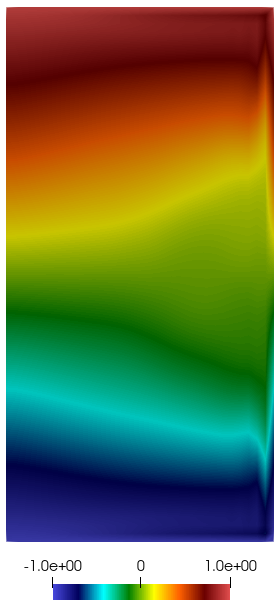}
    \put(6,92){\textcolor{white}{\scriptsize{2QGE, $32 \times 64$}}}
    \end{overpic}
    \begin{overpic}[percent,width=0.162\textwidth,grid=false]{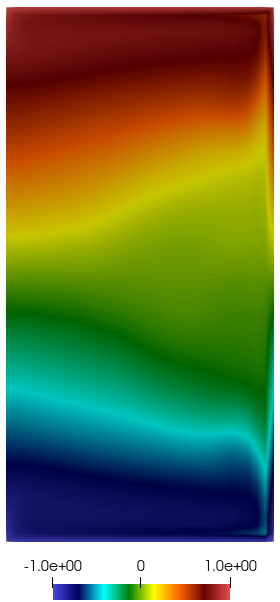}
    \put(5,92){\textcolor{white}{\scriptsize{2QG-$\alpha$, $32 \times 64$}}}
    \end{overpic}
        \begin{overpic}[percent,width=0.162\textwidth,grid=false]{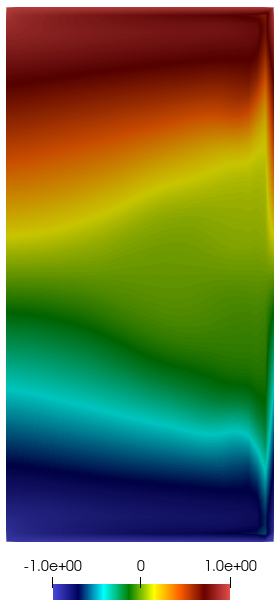}
    \put(1,92){\textcolor{white}{\scriptsize{2QG-NL-$\alpha$, $32 \times 64$}}}
    \end{overpic}\\
\begin{overpic}[percent,width=0.16\textwidth,grid=false]{img/mf-d0.5/DNS-q2mean0.5.png}
    \put(5,92){\textcolor{white}{\scriptsize{DNS, $256 \times 512$}}}
    \end{overpic}
    \begin{overpic}[percent,width=0.162\textwidth,grid=false]{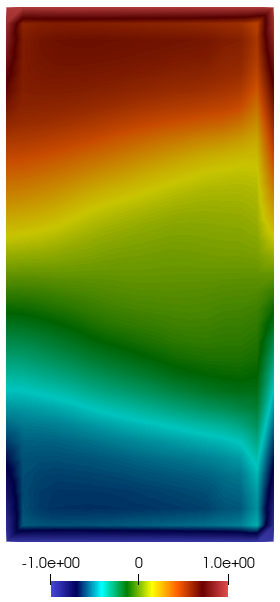}
    \put(6,92){\textcolor{white}{\scriptsize{2QGE, $16 \times 32$}}}
    \end{overpic}
    \begin{overpic}[percent,width=0.162\textwidth,grid=false]{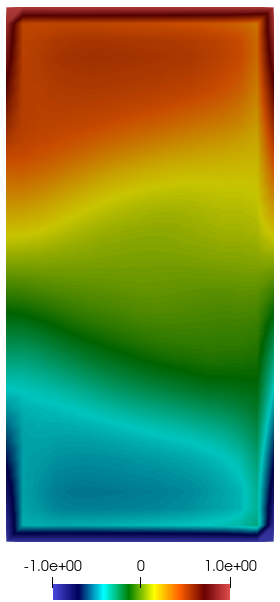}
    \put(5,92){\textcolor{white}{\scriptsize{2QG-$\alpha$, $16 \times 32$}}}
    \end{overpic}
        \begin{overpic}[percent,width=0.162\textwidth,grid=false]{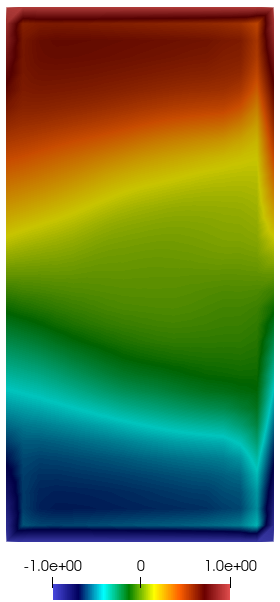}
    \put(1,92){\textcolor{white}{\scriptsize{2QG-NL-$\alpha$, $16 \times 32$}}}
       \end{overpic}\\
   \begin{overpic}[percent,width=0.16\textwidth,grid=false]{img/mf-d0.5/DNS-q2mean0.5.png}
    \put(5,92){\textcolor{white}{\scriptsize{DNS, $256 \times 512$}}}
    \end{overpic}
        \begin{overpic}[percent,width=0.162\textwidth,grid=false]{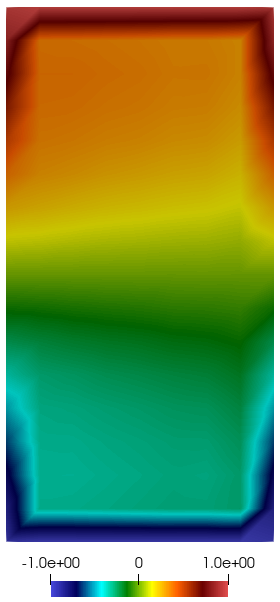}
    \put(6,92){\textcolor{white}{\scriptsize{2QGE, $8 \times 16$}}}
    \end{overpic}
    \begin{overpic}[percent,width=0.162\textwidth,grid=false]{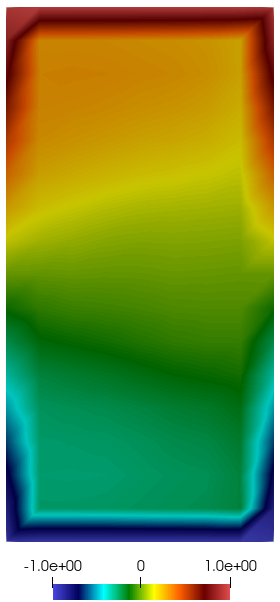}
    \put(5,92){\textcolor{white}{\scriptsize{2QG-$\alpha$, $8 \times 16$}}}
    \end{overpic}
    \begin{overpic}[percent,width=0.162\textwidth,grid=false]{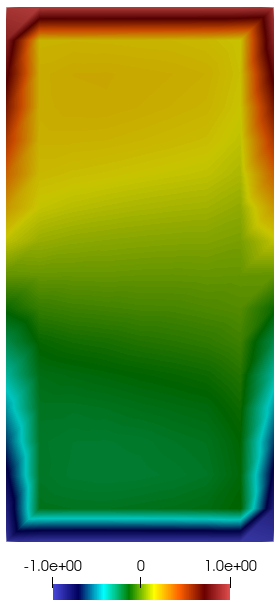}
    \put(2,92){\textcolor{white}{\scriptsize{2QG-NL-$\alpha$, $8 \times 16$}}}
       \end{overpic}
\caption{{Case 1:} Time-averaged potential vorticity of the bottom layer $\widetilde{q}_2$ computed by DNS (first column) and 2QGE with no filtering (second column), 2QG-$\alpha$ (third column), and 2QG-NL-$\alpha$ (fourth column) with different coarse meshes. Note that the color bar may differ
from one panel to the other.}
\label{fig:q2-case1-filtering}
\end{figure}

The time-averaged variables obtained with DNS, i.e., the 2QGE model with mesh $256 \times 512$, are shown on the first columns of Fig.~\ref{fig:psi1-case1-filtering}-\ref{fig:q2-case1-filtering}. 
We report the DNS solution on every row to ease the comparison with the solutions on coarser meshes.
We see that while $\widetilde{q}_1$ and $\widetilde{q}_2$ (first column of Fig.~\ref{fig:q1-case1-filtering} and \ref{fig:q2-case1-filtering}, respectively) present little interesting features, $\widetilde{\psi}_1$ (first column of Fig.~\ref{fig:psi1-case1-filtering}) presents two gyres with interesting shapes 
and $\widetilde{\psi}_2$ (first column of Fig.~\ref{fig:psi2-case1-filtering}) shows two well-defined inner gyres and two more faint outer gyres.

If one uses the 2QGE model with coarser meshes 
$32 \times 64$, $16 \times 32$, and $8 \times 16$ these  large-scale structures are reconstructed with the wrong shape, wrong magnitude or possibly disappear:
see the second columns in Fig.~\ref{fig:psi1-case1-filtering}-\ref{fig:psi2-case1-filtering}. 
With these coarse meshes, $\widetilde{q}_1$ and $\widetilde{q}_2$ are also not accurately reconstructed, however the magnitude is correctly captured and the differences are less evident, especially with mesh 
$32 \times 64$.
See the second columns in Fig.~\ref{fig:q1-case1-filtering}-\ref{fig:q2-case1-filtering}. Note that these coarser mesh all have a mesh size larger than $\delta_M$, thus it is expected that the
solver for the 2QGE (no filtering) yields
an inaccurate solution.

Next, we show that the time-averaged approximations obtained
with the same coarse meshes with the 2QG-$\alpha$ and 2QG-NL-$\alpha$ models ($\alpha = \sqrt{2} h$) are 
an improvement over the solutions given by the 2QGE model. Let us start from the stream functions.
For both the top and bottom layers, 
the 2QG-$\alpha$ produces a time-averaged stream function field with the correct number of gyres
and a range of values closer to that of the DNS. In particular, we note that the 
2QG-$\alpha$ model with mesh $32 \times 64$ provides rather accurate $\widetilde{\psi}_1$ 
and $\widetilde{\psi}_2$, especially towards the Eastern boundary. This accuracy degenerates
with the coarser meshes in terms of both gyre shape and value range. 
The 2QG-NL-$\alpha$ model produces time-averaged $\widetilde{\psi}_1$ fields that are a remarkable
improvement over the $\widetilde{\psi}_1$ given by the 2QG-$\alpha$ model. For all the coarser
meshes, the shape of the gyres gets closer to the shape observed in the DNS and, with the
exception of mesh $8 \times 16$, the value ranges matches well with the range in the DNS.
Gyres shape and value range for $\widetilde{\psi}_2$ given by the 2QG-NL-$\alpha$ model
with meshes $32 \times 64$ and $16 \times 32$ are also an improvement over the results
given by the 2QG-$\alpha$ model. However, that is not the case for mesh $8 \times 16$. 
We conclude that mesh $8 \times 16$ is too coarse to reconstruct the solutions features.
Hence, we will disregard it for case 2. Now, let us look at the vorticities in Fig.~\ref{fig:q1-case1-filtering}, \ref{fig:q2-case1-filtering}, third and fourth columns. 
The 2QG-$\alpha$ and 2QG-NL-$\alpha$ models improve the approximation of $\widetilde{q}_1$ near the Eastern boundary. There, in the top layer the flow coalesces, generating
a sharp transition layer between the red region (large positive values) and the 
blue region (large negative values). 
We see that both the 2QG-$\alpha$ and 2QG-NL-$\alpha$ models improve
the reconstruction of this transition layer, which is hard to capture with
coarse meshes. In particular, 
with meshes  $32 \times 64$ and $16 \times 32$ the 2QG-NL-$\alpha$ model 
gives a larger cyan (small negative values)
region in $\widetilde{q}_1$ than the 2QG-$\alpha$ model, which compares
favorably with the DNS. In $\widetilde{q}_2$ (see Fig.~\ref{fig:q2-case1-filtering}), 
there is a sharp boundary layer at the Easter boundary, 
where it has to satisfy the given boundary condition.
This boundary layer becomes hard to capture with coarse meshes. Nonetheless, the approximations of $\widetilde{q}_1$ computed with the 2QG-NL-$\alpha$ model and meshes $32 \times 64$ and $16 \times 32$
are pretty accurate. 

\begin{figure}[htb!]
\centering
    \begin{subfigure}{0.48\textwidth}
         \centering
         \includegraphics[width=\textwidth]{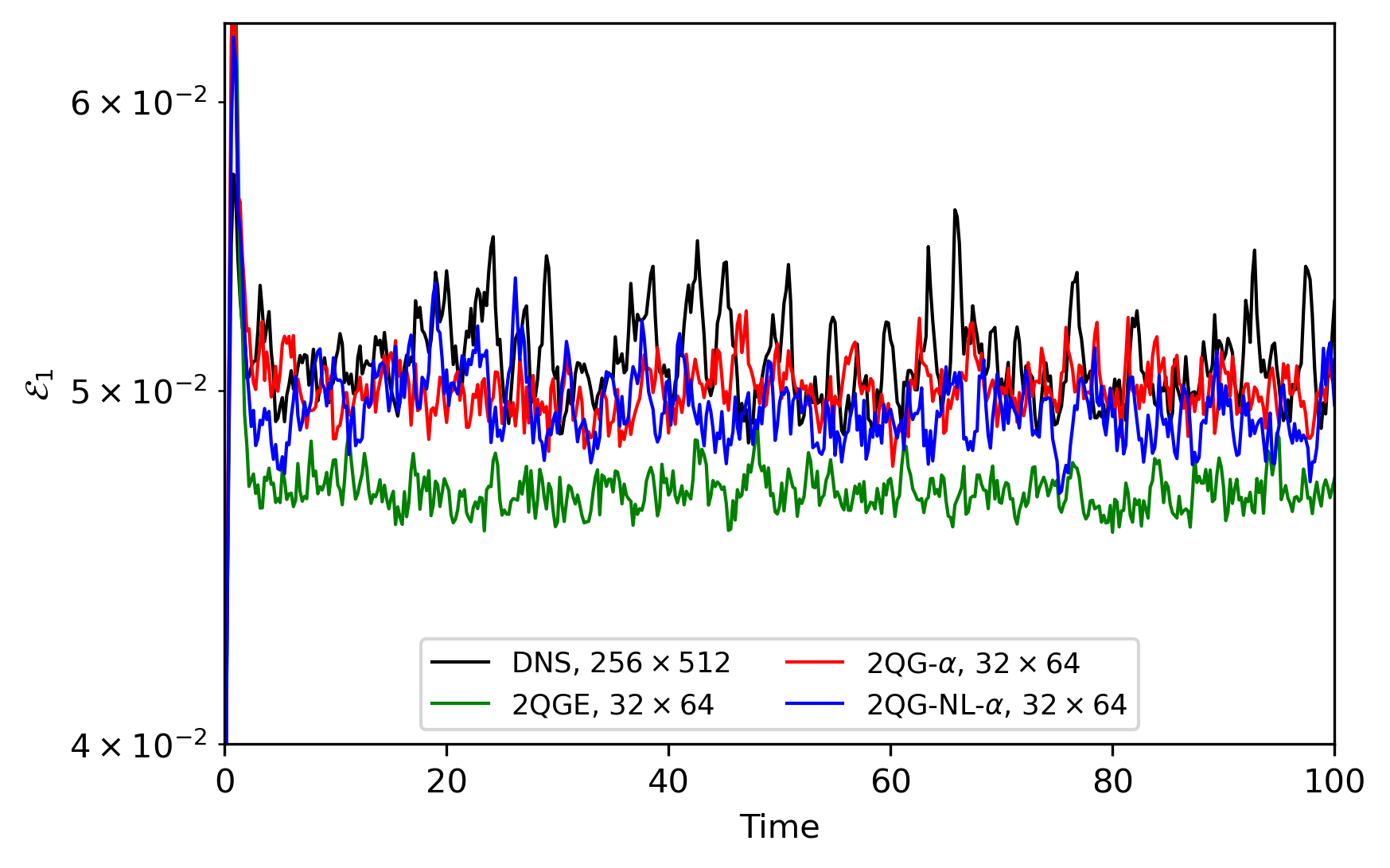}
     \end{subfigure}
    \begin{subfigure}{0.48\textwidth}
         \centering
         \includegraphics[width=\textwidth]{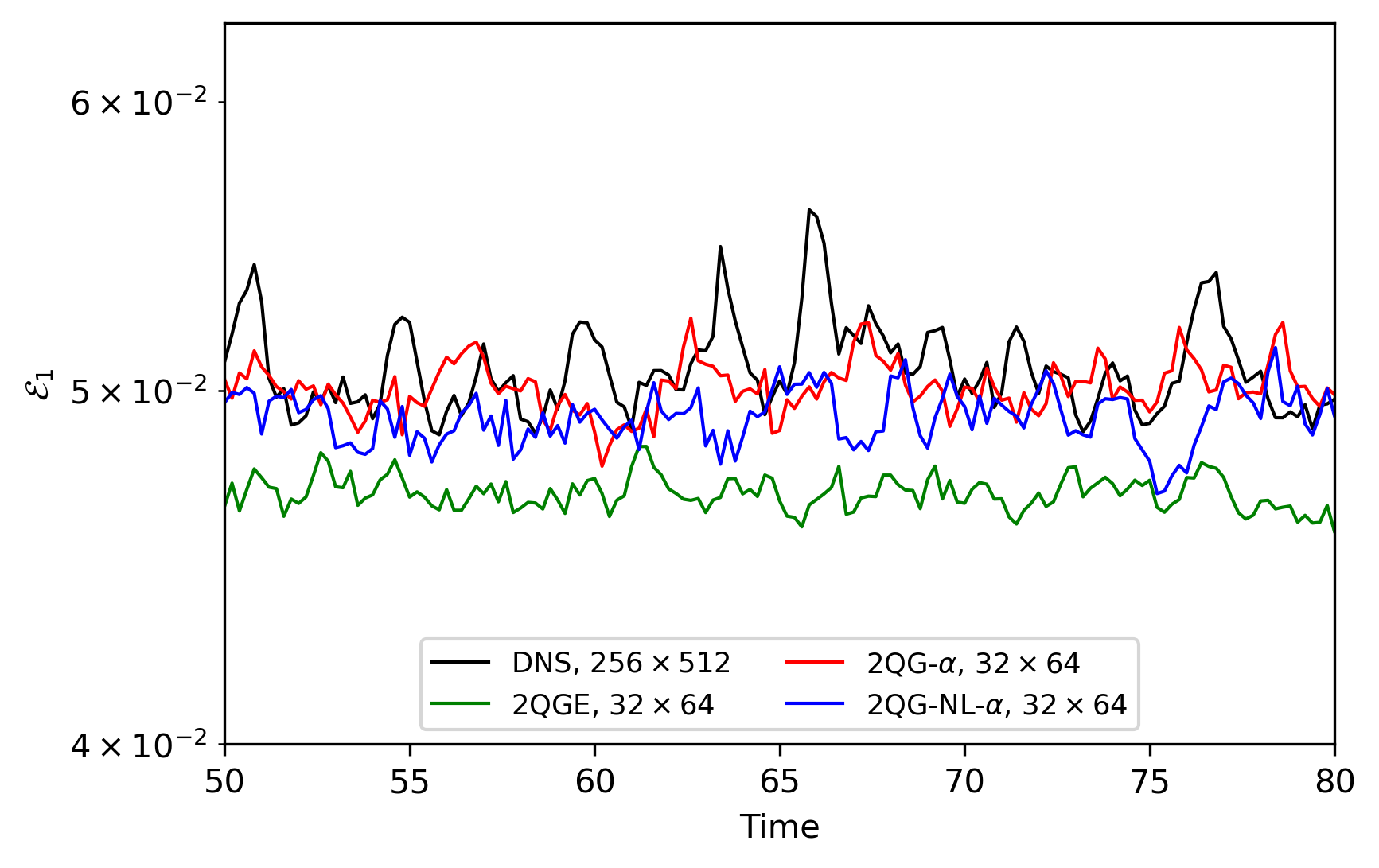}
     \end{subfigure}
    \begin{subfigure}{0.48\textwidth}
         \centering
         \includegraphics[width=\textwidth]{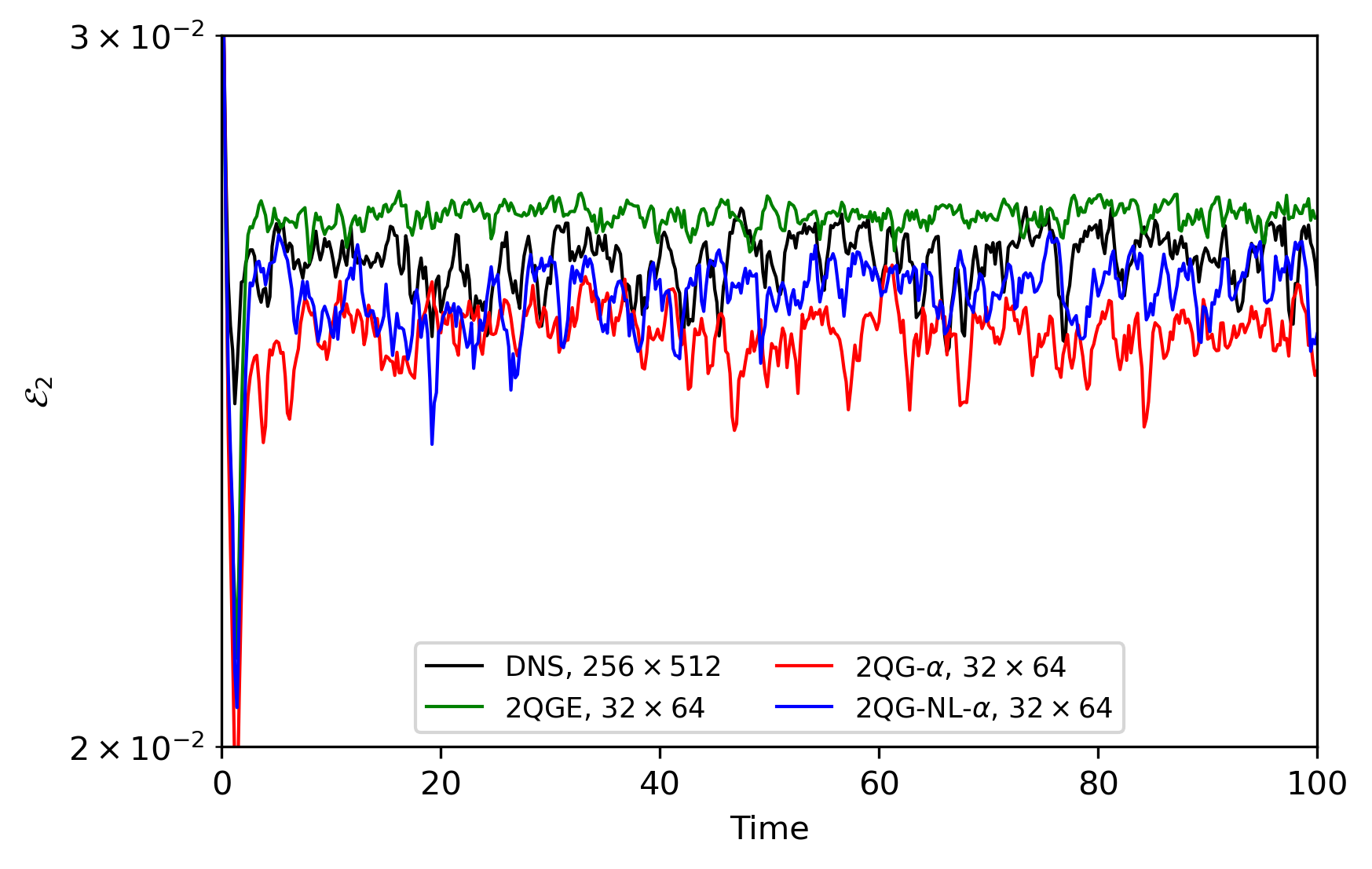}
     \end{subfigure}
    \begin{subfigure}{0.48\textwidth}
         \centering
         \includegraphics[width=\textwidth]{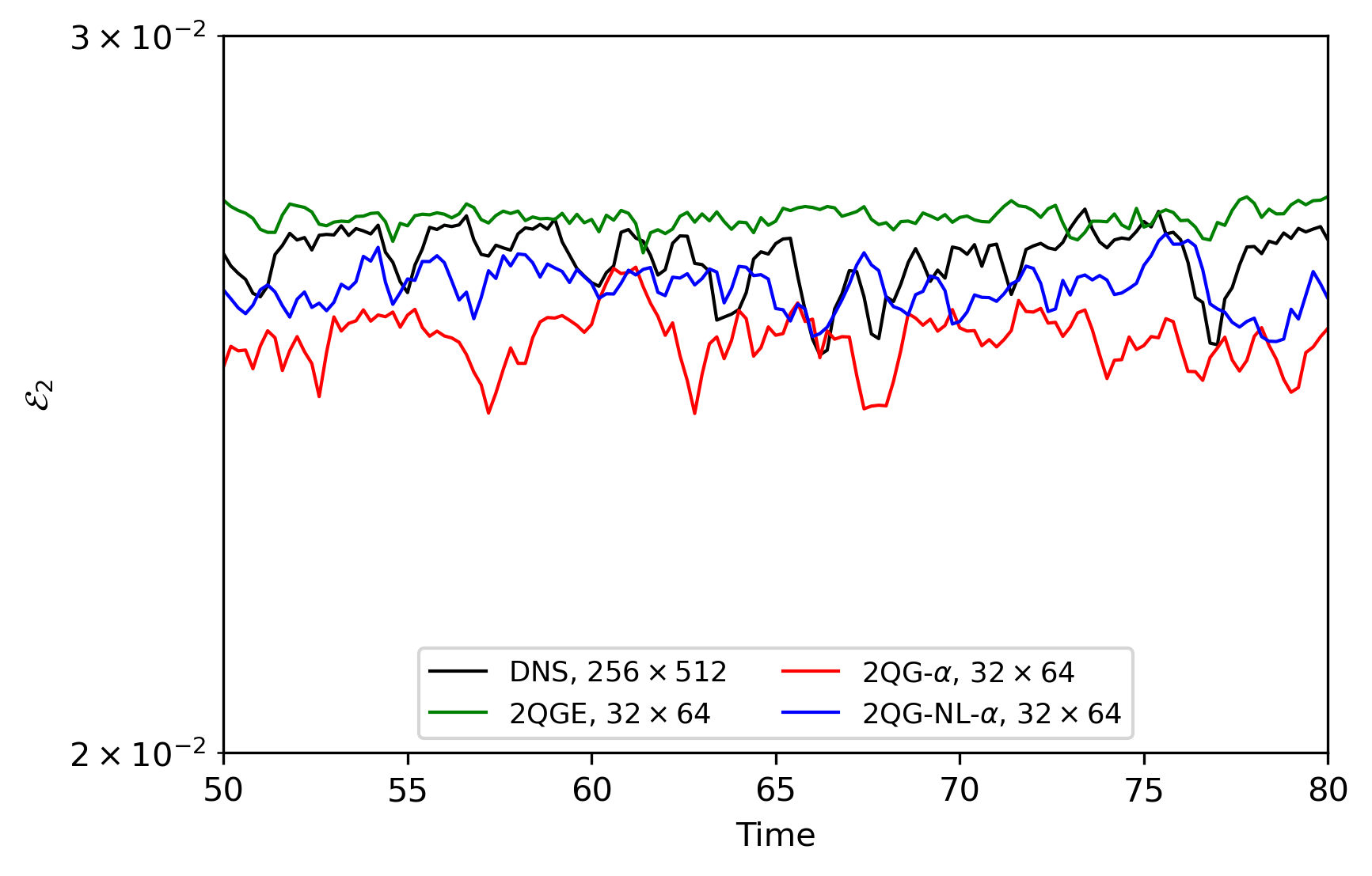}
     \end{subfigure}
\caption{Case 1: Time evolution of the enstrophy of the top layer $\mathcal{E}_1$ (top left) 
and the enstrophy of the bottom layer $\mathcal{E}_2$ (bottom left) 
for the DNS (mesh $256\times 512$), and 
2QGE model, 2QG-$\alpha$ model and 
2QG-NL-$\alpha$ model with coarse mesh $32 \times 64$. 
The panels on the right are zoomed-in views
of the panels on the left.}
\label{fig:case1-ens}
\end{figure}

Fig.~\ref{fig:case1-ens} shows the computed enstrophy \eqref{eq:entrophy} of the two layers 
for  the DNS, and the 2QGE, 2QG-$\alpha$, 2QG-NL-$\alpha$ models with a coarse mesh 
($32 \times 64$).
Notice that that the enstrophy computed with the 2QGE and mesh $32 \times 64$ consistently 
underestimates (resp., overestimates) the enstrophy computed with the DNS
in the top (resp., bottom) layer. With the same coarse mesh, both the 2QG-$\alpha$ and 
the 2QG-NL-$\alpha$ model improve $\mathcal{E}_1$ with respect to the 2QGE model. 
See top panels in Fig.~\ref{fig:case1-ens}.
However, it is only with the 2QG-NL-$\alpha$ model that one can improve $\mathcal{E}_2$.
In fact, $\mathcal{E}_2$ computed with the 2QG-$\alpha$ model tends to underestimate
the enstrophy of the bottom layer. See bottom panels in Fig.~\ref{fig:case1-ens}. A confirmation that the 2QG-NL-$\alpha$ model gives the best approximation
of the enstrophies computed by the DNS is given in Table \ref{tab:exp1-ens-32x64},
which reports the extrema of $\mathcal{E}_1$ and $\mathcal{E}_2$
for all the simulations in Fig.~\ref{fig:case1-ens}, together with the $L^2$ error. 

\begin{table}[h!]
    \centering
    \begin{tabular}{|c|c|c|c|c|c|c|}
       \hline
       \multirow{2}{*}{}  & \multirow{2}{*}{} & \multirow{2}{*}{DNS} & \multicolumn{3}{c|}{$32\times64$} \\ \cline{4-6} 
       \multirow{2}{*}{} & \multirow{2}{*}{} &  & 2QGE & 2QG-$\alpha$ & 2QG-NL-$\alpha$ \\ \hline
       \multirow{3}{*}{$\mathcal{E}_1$} & average & 5.094E-02 & 4.704E-02 & {5.011E-02} & 4.950E-02 \\ \cline{2-6}
       & max & 5.732E-02 & 6.736E-02 & 6.544E-02 & {6.249E-02} \\ \cline{2-6}
       & $L^2$ error & - & 9.616E-02 & {4.343E-02} & 4.930E-02 \\ \hline
       \multirow{3}{*}{$\mathcal{E}_2$} & average & 2.640E-02  & 2.703E-02 & 2.526E-02 & {2.587E-02} \\ \cline{2-6}
       & min & 2.432E-02 & {2.103E-02} & 1.917E-02 & 2.045E-02 \\ \cline{2-6}
       & $L^2$ error & - & 1.950E-02 & 2.999E-02 & {1.831E-02} \\ \hline
    \end{tabular}
    \caption{Case 1: extrema of $\mathcal{E}_1$ and $\mathcal{E}_2$ for the
    DNS (mesh $256\times512$) and the 2QGE, 2QG-$\alpha$, 2QG-NL-$\alpha$ models 
    with coarse mesh $32\times64$. For the 2QGE, 2QG-$\alpha$, 2QG-NL-$\alpha$ models, we also
    report the associated $L^2$ errors.}
    \label{tab:exp1-ens-32x64}
\end{table}

We conclude this subsection with a comment on
the computational cost of the simulations presented in this section. All simulations were run on a computing server equipped with Intel\textsuperscript{\textregistered} Xeon\textsuperscript{\textregistered} processor and 128GB of RAM, with a 64-bit version of Linux.
Table \ref{tab:cpu_time} reports all the computational costs
and the speed up allowed by the use of filtering with coarser meshes. We do not report the speed up for the 2QGE model
with coarser meshes since we have seen in Fig.~\ref{fig:psi1-case1-filtering}-\ref{fig:q2-case1-filtering} that the corresponding solutions do not compare well with the DNS. The first thing to notice
in Table \ref{tab:cpu_time} is that the DNS takes more than 8 days to complete. When using the 2QG-$\alpha$ model,
the computational cost drops to about 1 hour and 40 minutes
(118 times faster than the DNS) 
with mesh $32 \times 64$ and about 39 minutes 
(299 times faster than the DNS)
with mesh $16 \times 32$.
We recall that the 2QG-$\alpha$ model requires the solution of two additional 
linear filter problems. The computational cost for those two
additional problems adds up to only 20 minutes for
mesh $32 \times 64$ and 8 minutes for mesh $16 \times 32$. 
Compare the times in Table \ref{tab:cpu_time} for the 2QGE model 
and the 2QG-$\alpha$ model for a given mesh. 
Fig.~\ref{fig:psi1-case1-filtering}-\ref{fig:q2-case1-filtering}
clearly indicated that the solutions given by the 2QG-NL-$\alpha$ model 
are more accurate than the solutions provided by the 2QG-$\alpha$ model, although one has to pay the price of solving two
nonlinear filter problems instead of two linear filter problems.
Table \ref{tab:cpu_time} shows that this price is just
13 minutes for mesh $32 \times 64$ and 5 minutes for mesh
$16 \times 32$.

\begin{table}[htb!]
    \centering
    \begin{tabular}{lllc}
        \hline 
        Model & Mesh & CPU Time & Speed Up Factor \\
        \hline
        DNS & $256\times512$ & 8d 4h 41m & - \\ 
        2QGE & $32\times64$ & 1h 20m & - \\ 
        2QG-$\alpha$ & $32\times64$ & 1h 40m & 118 \\ 
        2QG-NL-$\alpha$ & $32\times64$ & 1h 53m & 104 \\ 
        2QGE & $16\times32$ & 31m 40s & - \\ 
        2QG-$\alpha$ & $16\times32$ & 39m 30s & 299 \\ 
        2QG-NL-$\alpha$ & $16\times32$ & 44m 13s & 268 \\ 
        \hline
    \end{tabular}
    \caption{Case 1: computational time required by the DNS and
    the 2QGE, 2QG-$\alpha$, and 2QG-NL-$\alpha$ models with coarser meshes
    $32\times64$ and $16\times32$. For the algorithms with filtering, the
    speed up factor with respect to the DNS is also reported.}
    \label{tab:cpu_time}
\end{table}

\subsubsection{Numerical results for case 2} \label{sec:experiment2}

Case 2 has the same set-up as case 1, and in particular
same diffusivity, but larger $\sigma$ (more friction with the bottom of the ocean than in case 1) and smaller $\delta$ (thinner top layer than in case 1).
Recall that smaller $\delta$ means that the dynamics
in the bottom layer affects more strongly the dynamics
of the top layer, while the dynamics
in the bottom layer is less affected by the top layer.
We will proceed in the same manner as for case 1, i.e., we take
the DNS to be the true solution and consider three coarse 
meshes to see if, through the use of filtering, we can obtain
a solution close to the true one. 

\begin{figure}
\centering
    \begin{overpic}[percent,width=0.16\textwidth,grid=false]{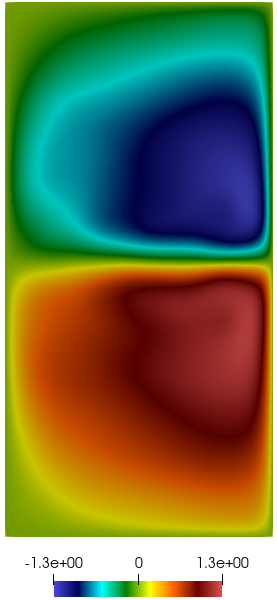}
    \put(5,92){\textcolor{white}{\scriptsize{DNS, $256 \times 512$}}}
    \end{overpic}
    \begin{overpic}[percent,width=0.16\textwidth,grid=false]{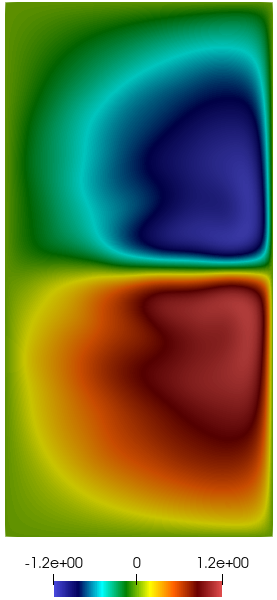}
    \put(6,92){\textcolor{white}{\scriptsize{2QGE, $64 \times 128$}}}
    \end{overpic}
    \begin{overpic}[percent,width=0.16\textwidth,grid=false]{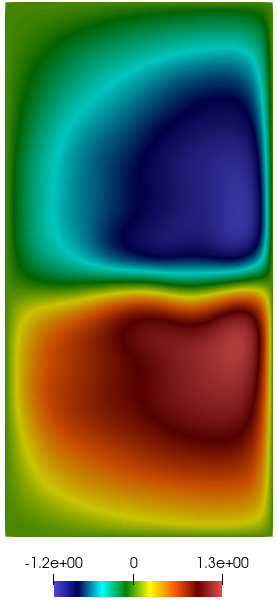}
    \put(5,92){\textcolor{white}{\scriptsize{2QG-$\alpha$, $64 \times 128$}}}
    \end{overpic}
        \begin{overpic}[percent,width=0.16\textwidth,grid=false]{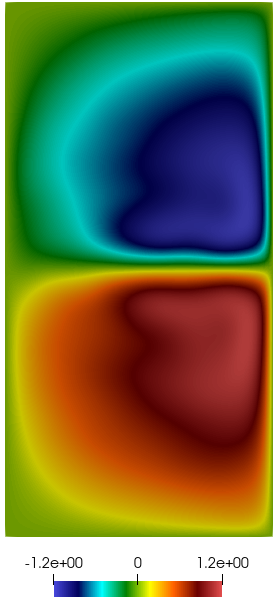}
    \put(2,92){\textcolor{white}{\tiny{2QG-NL-$\alpha$, $64 \times 128$}}}
    \end{overpic}\\
    \begin{overpic}[percent,width=0.16\textwidth,grid=false]{img/mf-case2/DNS-psi1mean2.png}
    \put(5,92){\textcolor{white}{\scriptsize{DNS, $256 \times 512$}}}
    \end{overpic}
    \begin{overpic}[percent,width=0.16\textwidth,grid=false]{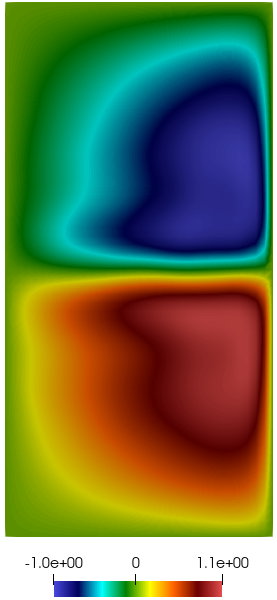}
    \put(6,92){\textcolor{white}{\scriptsize{2QGE, $32 \times 64$}}}
    \end{overpic}
    \begin{overpic}[percent,width=0.16\textwidth,grid=false]{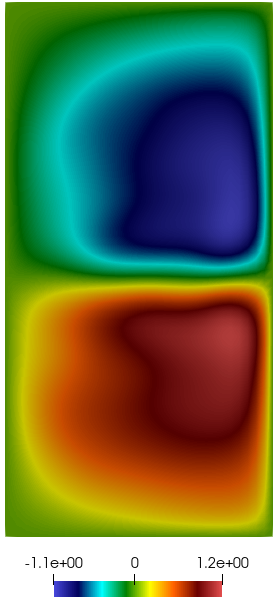}
    \put(5,92){\textcolor{white}{\scriptsize{2QG-$\alpha$, $32 \times 64$}}}
    \end{overpic}
        \begin{overpic}[percent,width=0.16\textwidth,grid=false]{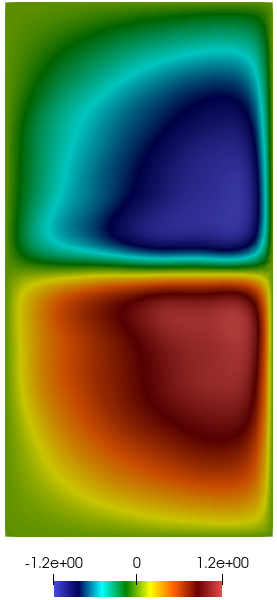}
    \put(1,92){\textcolor{white}{\scriptsize{2QG-NL-$\alpha$, $32 \times 64$}}}
    \end{overpic}\\
    \begin{overpic}[percent,width=0.16\textwidth,grid=false]{img/mf-case2/DNS-psi1mean2.png}
    \put(5,92){\textcolor{white}{\scriptsize{DNS, $256 \times 512$}}}
    \end{overpic}
    \begin{overpic}[percent,width=0.16\textwidth,grid=false]{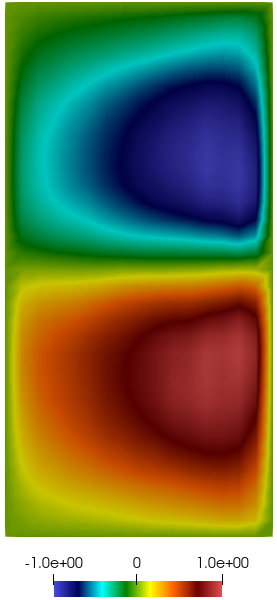}
    \put(6,92){\textcolor{white}{\scriptsize{2QGE, $16 \times 32$}}}
    \end{overpic}
    \begin{overpic}[percent,width=0.16\textwidth,grid=false]{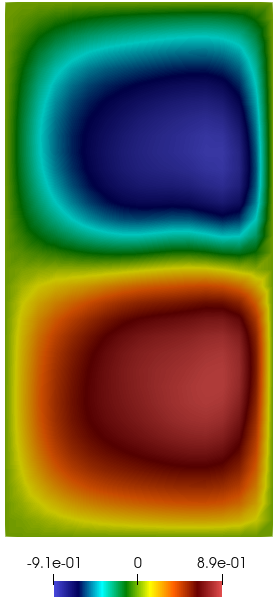}
    \put(5,92){\textcolor{white}{\scriptsize{2QG-$\alpha$, $16 \times 32$}}}
    \end{overpic}
        \begin{overpic}[percent,width=0.16\textwidth,grid=false]{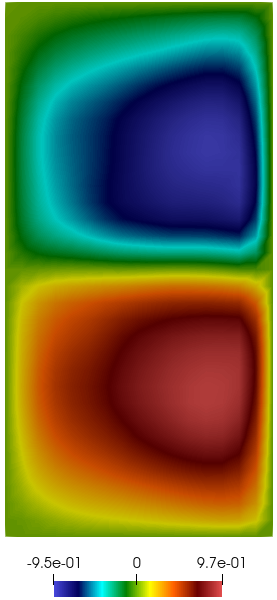}
    \put(1,92){\textcolor{white}{\scriptsize{2QG-NL-$\alpha$, $16 \times 32$}}}
    \end{overpic}\\
\caption{{Case 2:} Time-averaged stream function of the top layer $\widetilde{\psi}_1$ computed by DNS (first column) and 2QGE with no filtering (second column), 2QG-$\alpha$ (third column), and 2QG-NL-$\alpha$ (fourth column) with different coarse meshes. Note that the color bar may differ
from one panel to the other.}
\label{fig:psi1-case2-filtering}
\end{figure}

\begin{figure}
\centering
    \begin{overpic}[percent,width=0.16\textwidth,grid=false]{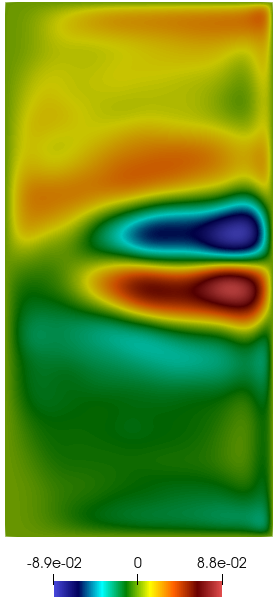}
    \put(5,92){\textcolor{white}{\scriptsize{DNS, $256 \times 512$}}}
    \end{overpic}
    \begin{overpic}[percent,width=0.16\textwidth,grid=false]{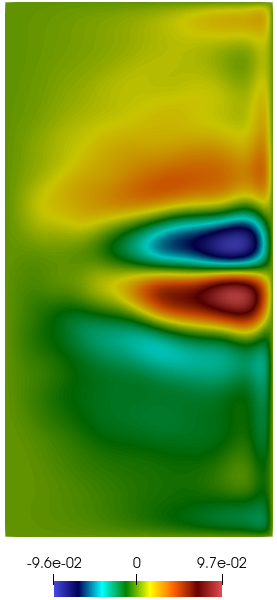}
    \put(6,92){\textcolor{white}{\scriptsize{2QGE, $64 \times 128$}}}
    \end{overpic}
    \begin{overpic}[percent,width=0.16\textwidth,grid=false]{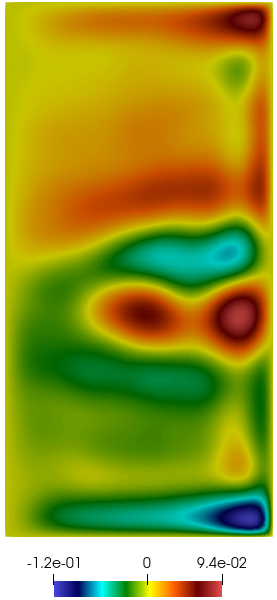}
    \put(5,92){\textcolor{white}{\scriptsize{2QG-$\alpha$, $64 \times 128$}}}
    \end{overpic}
        \begin{overpic}[percent,width=0.16\textwidth,grid=false]{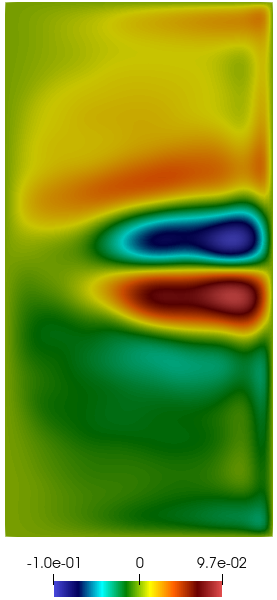}
    \put(2,92){\textcolor{white}{\tiny{2QG-NL-$\alpha$, $64 \times 128$}}}
    \end{overpic}\\
    \begin{overpic}[percent,width=0.16\textwidth,grid=false]{img/mf-case2/DNS-psi2mean2.png}
    \put(5,92){\textcolor{white}{\scriptsize{DNS, $256 \times 512$}}}
    \end{overpic}
    \begin{overpic}[percent,width=0.16\textwidth,grid=false]{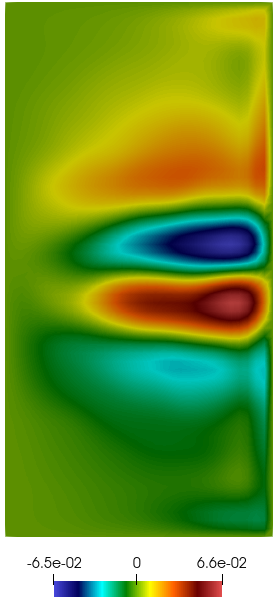}
    \put(6,92){\textcolor{white}{\scriptsize{2QGE, $32 \times 64$}}}
    \end{overpic}
    \begin{overpic}[percent,width=0.16\textwidth,grid=false]{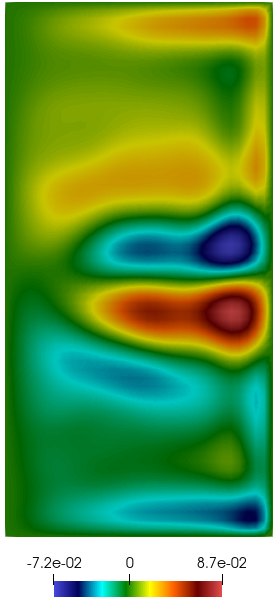}
    \put(5,92){\textcolor{white}{\scriptsize{2QG-$\alpha$, $32 \times 64$}}}
    \end{overpic}
        \begin{overpic}[percent,width=0.16\textwidth,grid=false]{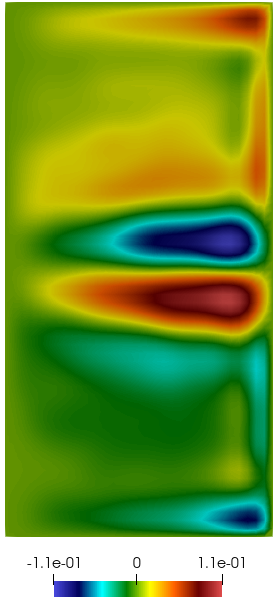}
    \put(1,92){\textcolor{white}{\scriptsize{2QG-NL-$\alpha$, $32 \times 64$}}}
    \end{overpic}\\
    \begin{overpic}[percent,width=0.16\textwidth,grid=false]{img/mf-case2/DNS-psi2mean2.png}
    \put(5,92){\textcolor{white}{\scriptsize{DNS, $256 \times 512$}}}
    \end{overpic}
    \begin{overpic}[percent,width=0.16\textwidth,grid=false]{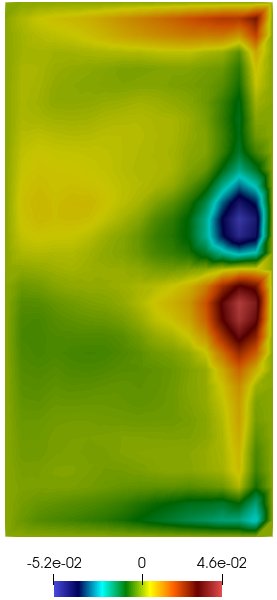}
    \put(6,92){\textcolor{white}{\scriptsize{2QGE, $16 \times 32$}}}
    \end{overpic}
    \begin{overpic}[percent,width=0.16\textwidth,grid=false]{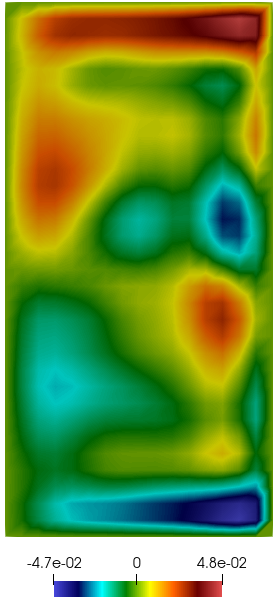}
    \put(5,92){\textcolor{white}{\scriptsize{2QG-$\alpha$, $16 \times 32$}}}
    \end{overpic}
        \begin{overpic}[percent,width=0.16\textwidth,grid=false]{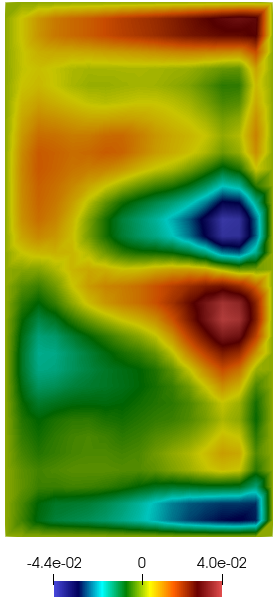}
    \put(1,92){\textcolor{white}{\scriptsize{2QG-NL-$\alpha$, $16 \times 32$}}}
    \end{overpic}\\
\caption{{Case 2:} Time-averaged stream function of the bottom layer $\widetilde{\psi}_2$ computed by DNS (first column) and 2QGE with no filtering (second column), 2QG-$\alpha$ (third column), and 2QG-NL-$\alpha$ (fourth column) with different coarse meshes. Note that the color bar may differ
from one panel to the other.}
\label{fig:psi2-case2-filtering}
\end{figure}

\begin{figure}
\centering
    \begin{overpic}[percent,width=0.16\textwidth,grid=false]{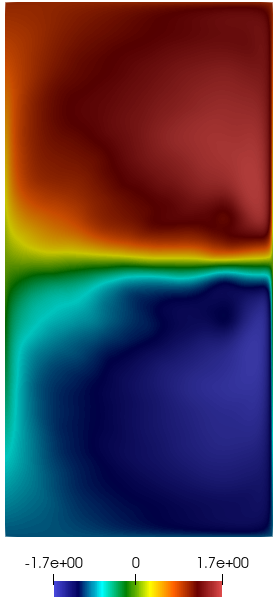}
    \put(5,92){\textcolor{white}{\scriptsize{DNS, $256 \times 512$}}}
    \end{overpic}
    \begin{overpic}[percent,width=0.16\textwidth,grid=false]{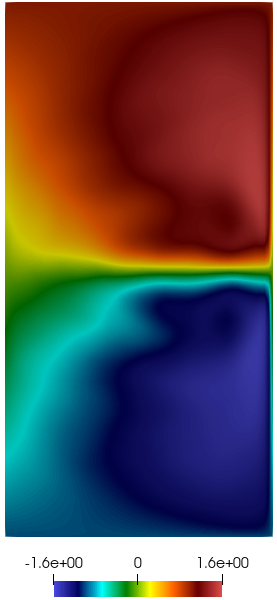}
    \put(6,92){\textcolor{white}{\scriptsize{2QGE, $64 \times 128$}}}
    \end{overpic}
    \begin{overpic}[percent,width=0.16\textwidth,grid=false]{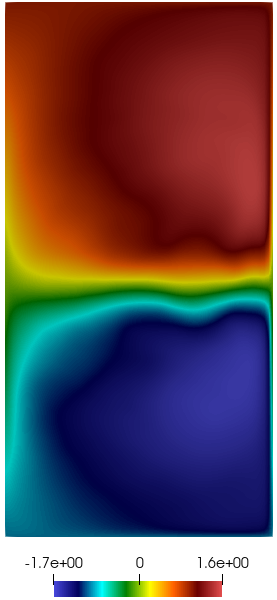}
    \put(5,92){\textcolor{white}{\scriptsize{2QG-$\alpha$, $64 \times 128$}}}
    \end{overpic}
        \begin{overpic}[percent,width=0.16\textwidth,grid=false]{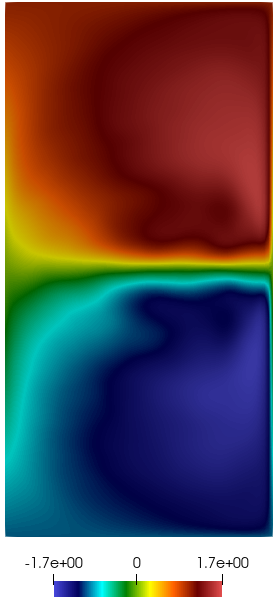}
    \put(2,92){\textcolor{white}{\tiny{2QG-NL-$\alpha$, $64 \times 128$}}}
    \end{overpic}\\
    \begin{overpic}[percent,width=0.16\textwidth,grid=false]{img/mf-case2/DNS-q1mean2.png}
    \put(5,92){\textcolor{white}{\scriptsize{DNS, $256 \times 512$}}}
    \end{overpic}
    \begin{overpic}[percent,width=0.16\textwidth,grid=false]{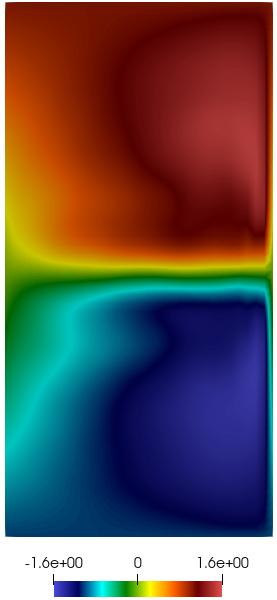}
    \put(6,92){\textcolor{white}{\scriptsize{2QGE, $32 \times 64$}}}
    \end{overpic}
    \begin{overpic}[percent,width=0.16\textwidth,grid=false]{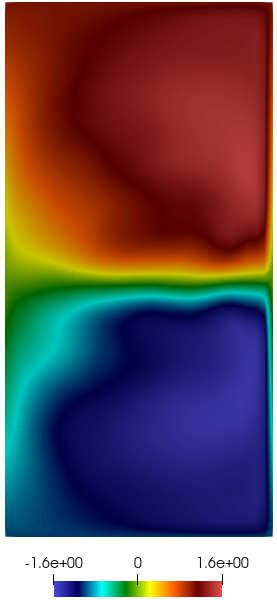}
    \put(5,92){\textcolor{white}{\scriptsize{2QG-$\alpha$, $32 \times 64$}}}
    \end{overpic}
        \begin{overpic}[percent,width=0.16\textwidth,grid=false]{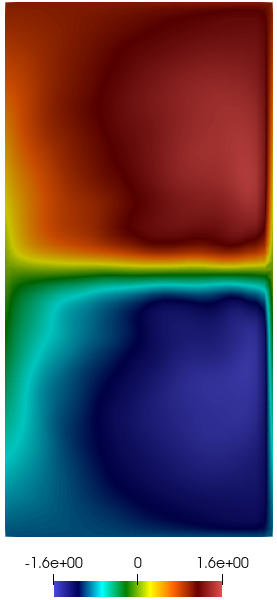}
    \put(1,92){\textcolor{white}{\scriptsize{2QG-NL-$\alpha$, $32 \times 64$}}}
    \end{overpic}\\
    \begin{overpic}[percent,width=0.16\textwidth,grid=false]{img/mf-case2/DNS-q1mean2.png}
    \put(5,92){\textcolor{white}{\scriptsize{DNS, $256 \times 512$}}}
    \end{overpic}
    \begin{overpic}[percent,width=0.16\textwidth,grid=false]{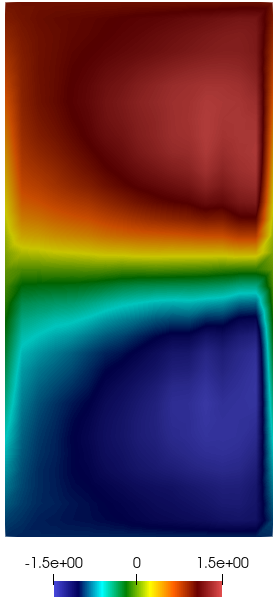}
    \put(6,92){\textcolor{white}{\scriptsize{2QGE, $16 \times 32$}}}
    \end{overpic}
    \begin{overpic}[percent,width=0.16\textwidth,grid=false]{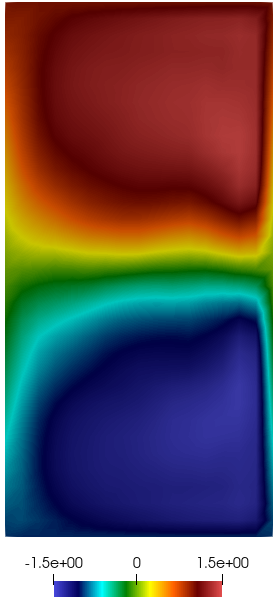}
    \put(5,92){\textcolor{white}{\scriptsize{2QG-$\alpha$, $16 \times 32$}}}
    \end{overpic}
        \begin{overpic}[percent,width=0.16\textwidth,grid=false]{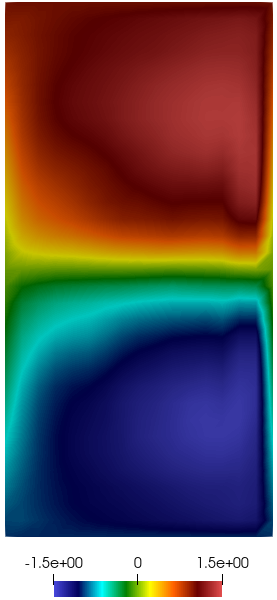}
    \put(1,92){\textcolor{white}{\scriptsize{2QG-NL-$\alpha$, $16 \times 32$}}}
    \end{overpic}\\
\caption{{Case 2:} Time-averaged potential vorticity of the top layer $\widetilde{q}_1$ computed by DNS (first column) and 2QGE with no filtering (second column), 2QG-$\alpha$ (third column), and 2QG-NL-$\alpha$ (fourth column) with different coarse meshes. Note that the color bar may differ
from one panel to the other.}
\label{fig:q1-case2-filtering}
\end{figure}

\begin{figure}
\centering
    \begin{overpic}[percent,width=0.16\textwidth,grid=false]{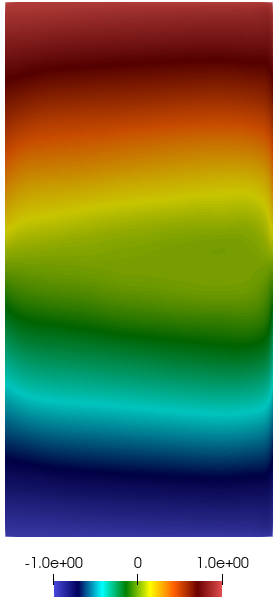}
    \put(5,92){\textcolor{white}{\scriptsize{DNS, $256 \times 512$}}}
    \end{overpic}
    \begin{overpic}[percent,width=0.16\textwidth,grid=false]{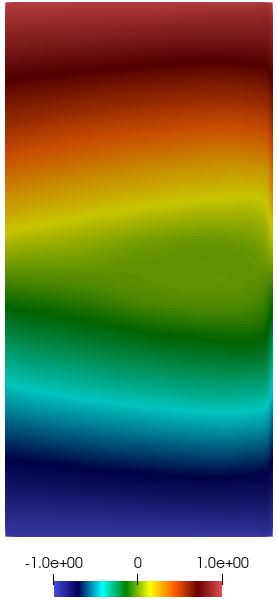}
    \put(6,92){\textcolor{white}{\scriptsize{2QGE, $64 \times 128$}}}
    \end{overpic}
    \begin{overpic}[percent,width=0.16\textwidth,grid=false]{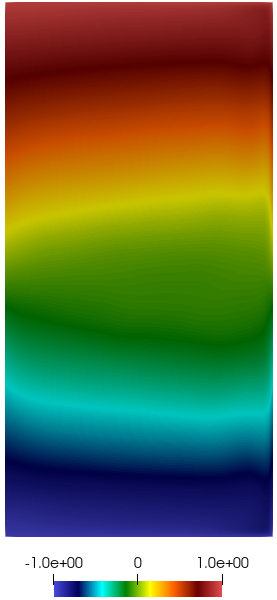}
    \put(5,92){\textcolor{white}{\scriptsize{2QG-$\alpha$, $64 \times 128$}}}
    \end{overpic}
        \begin{overpic}[percent,width=0.16\textwidth,grid=false]{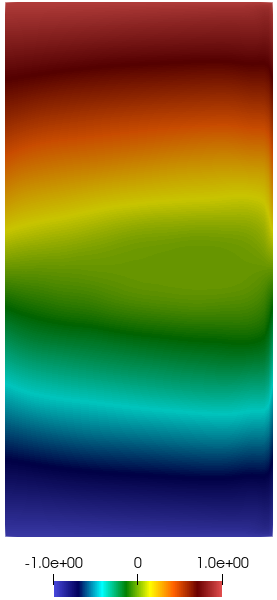}
    \put(2,92){\textcolor{white}{\tiny{2QG-NL-$\alpha$, $64 \times 128$}}}
    \end{overpic}\\
    \begin{overpic}[percent,width=0.16\textwidth,grid=false]{img/mf-case2/DNS-q2mean2.png}
    \put(5,92){\textcolor{white}{\scriptsize{DNS, $256 \times 512$}}}
    \end{overpic}
    \begin{overpic}[percent,width=0.16\textwidth,grid=false]{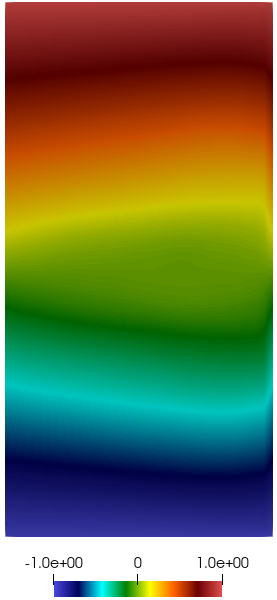}
    \put(6,92){\textcolor{white}{\scriptsize{2QGE, $32 \times 64$}}}
    \end{overpic}
    \begin{overpic}[percent,width=0.16\textwidth,grid=false]{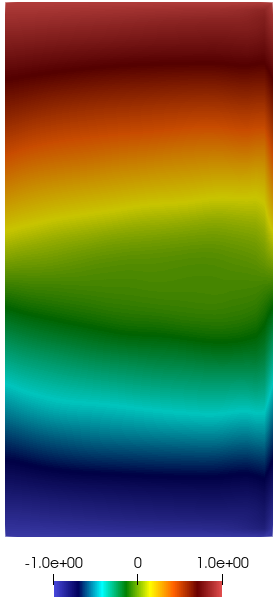}
    \put(5,92){\textcolor{white}{\scriptsize{2QG-$\alpha$, $32 \times 64$}}}
    \end{overpic}
        \begin{overpic}[percent,width=0.16\textwidth,grid=false]{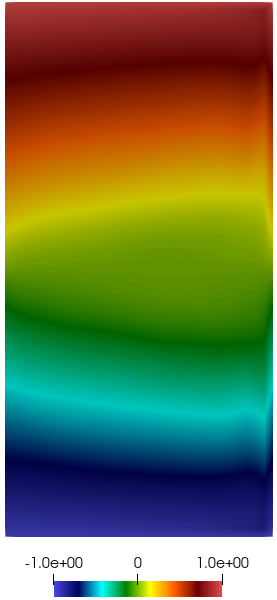}
    \put(1,92){\textcolor{white}{\scriptsize{2QG-NL-$\alpha$, $32 \times 64$}}}
    \end{overpic}\\
    \begin{overpic}[percent,width=0.16\textwidth,grid=false]{img/mf-case2/DNS-q2mean2.png}
    \put(5,92){\textcolor{white}{\scriptsize{DNS, $256 \times 512$}}}
    \end{overpic}
    \begin{overpic}[percent,width=0.16\textwidth,grid=false]{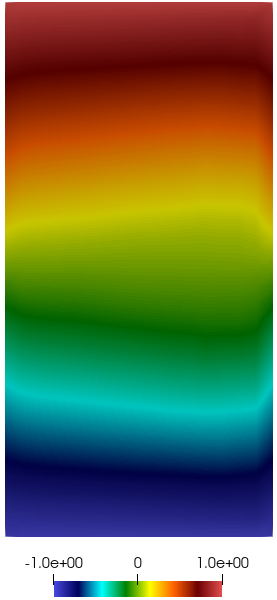}
    \put(6,92){\textcolor{white}{\scriptsize{2QGE, $16 \times 32$}}}
    \end{overpic}
    \begin{overpic}[percent,width=0.16\textwidth,grid=false]{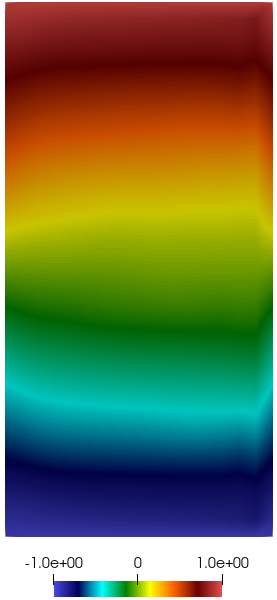}
    \put(5,92){\textcolor{white}{\scriptsize{2QG-$\alpha$, $16 \times 32$}}}
    \end{overpic}
        \begin{overpic}[percent,width=0.16\textwidth,grid=false]{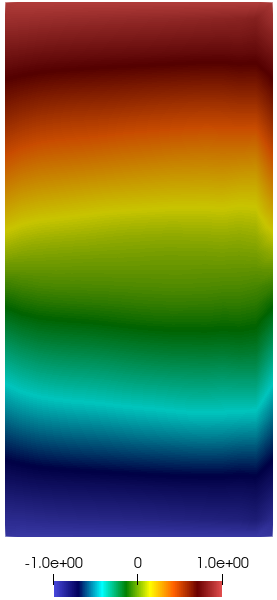}
    \put(1,92){\textcolor{white}{\scriptsize{2QG-NL-$\alpha$, $16 \times 32$}}}
    \end{overpic}\\
\caption{{Case 2:} Time-averaged potential vorticity of the bottom layer $\widetilde{q}_2$ computed by DNS (first column) and 2QGE with no filtering (second column), 2QG-$\alpha$ (third column), and 2QG-NL-$\alpha$ (fourth column) with different coarse meshes. Note that the color bar may differ
from one panel to the other.}
\label{fig:q2-case2-filtering}
\end{figure}

The first columns of Fig. \ref{fig:psi1-case2-filtering}-\ref{fig:q2-case2-filtering} display the time-averaged fields computed by the 
DNS, i.e., with the 2QGE model and mesh $256\times512$. 
We see that $\widetilde{\psi}_1$ shows two large gyres with a very
different shape with respect to case 1, while  $\widetilde{\psi}_2$ exhibits two smaller elongated inner gyres
and two larger outer gyres with a less defined shape. Concerning the time-averaged potential vorticities, note that the regions with
large vorticity in absolute value (red and blue) have expanded
with respect to case 1 in both the top and bottom layer, 
generating narrower transition regions (green)
at the center of the basin.
Given that the flow patterns are more challening in this
second case, we choose to consider meshes $16\times32$ and $32\times64$, 
which were used also for case 1, and one finer meshes, i.e., $64\times128$. 
While mesh $64\times128$ has a mesh size (=0.016) smaller than $\delta_M$ (=0.026), we will see that, even with this mesh, the solution given by the 2QGE model is not accurate.
The 2QGE model with the three coarser meshes
cannot correctly reconstruct
the shape of the gyres in $\widetilde{\psi}_1$ and $\widetilde{\psi}_2$ (second column in Fig.~\ref{fig:psi1-case2-filtering} and \ref{fig:psi2-case2-filtering}). Notice
that also the magnitude of $\widetilde{\psi}_1$ and $\widetilde{\psi}_2$ progressively degenerates
are the mesh gets coarser. Like for case 1, 
the 2QGE model with coarser meshes does a better job with $\widetilde{q}_1$ and $\widetilde{q}_2$. 
In fact, the second column in Fig.~\ref{fig:q2-case2-filtering} shows that $\widetilde{q}_2$ is very well reconstructed with all the three meshes. 
As for $\widetilde{q}_1$, we see that the green transition zone is rather well captured with 
all the three meshes, although with mesh $16\times 32$ it gets straightened out. 
See second column in Fig.~\ref{fig:q1-case2-filtering}.

Let us now have a look at the results given by 
the 2QG-$\alpha$ and 2QG-NL-$\alpha$ models  ($\alpha = h$) to 
see if and how they improve the results given by the 2QGE model for a given mesh. The first thing to notice is that the range of values  of $\widetilde{\psi}_1$ and $\widetilde{\psi}_2$ 
gets closer to the true one when using the 
2QG-$\alpha$ model. However, the shape of 
the gyres is not well reconstructed, especially as the mesh gets coarser. With the use of the nonlinear filter, i.e., the 2QG-NL-$\alpha$ model, 
these large scale structures in $\widetilde{\psi}_1$ and $\widetilde{\psi}_2$
are more accurately reconstructed, while further
improving the magnitude in most cases, the main exception being mesh $16 \times 32$. In fact, 
just like for case 1 we concluded that mesh
$8 \times 16$ was too coarse to capture
the time-averaged fields, we can conclude that mesh $16 \times 32$ is too coarse for case 2. 
In each case, there is a level of coarseness beyond which the 2QG-NL-$\alpha$ cannot improve the results.  From Fig.~\ref{fig:q1-case2-filtering}, we observe that
the 2QG-NL-$\alpha$ model provides the most accurate
$\widetilde{q}_1$ for mesh $64 \times 128$
and $32 \times 64$. Since $\widetilde{q}_2$ is
already well approximated by the 2QGE model, the action of both linear and nonlinear filters is limited. Compare the last three columns in Fig.~\ref{fig:q2-case2-filtering}.

\begin{figure}[htb!]
\centering
    \begin{subfigure}{0.48\textwidth}
         \centering
         \includegraphics[width=\textwidth]{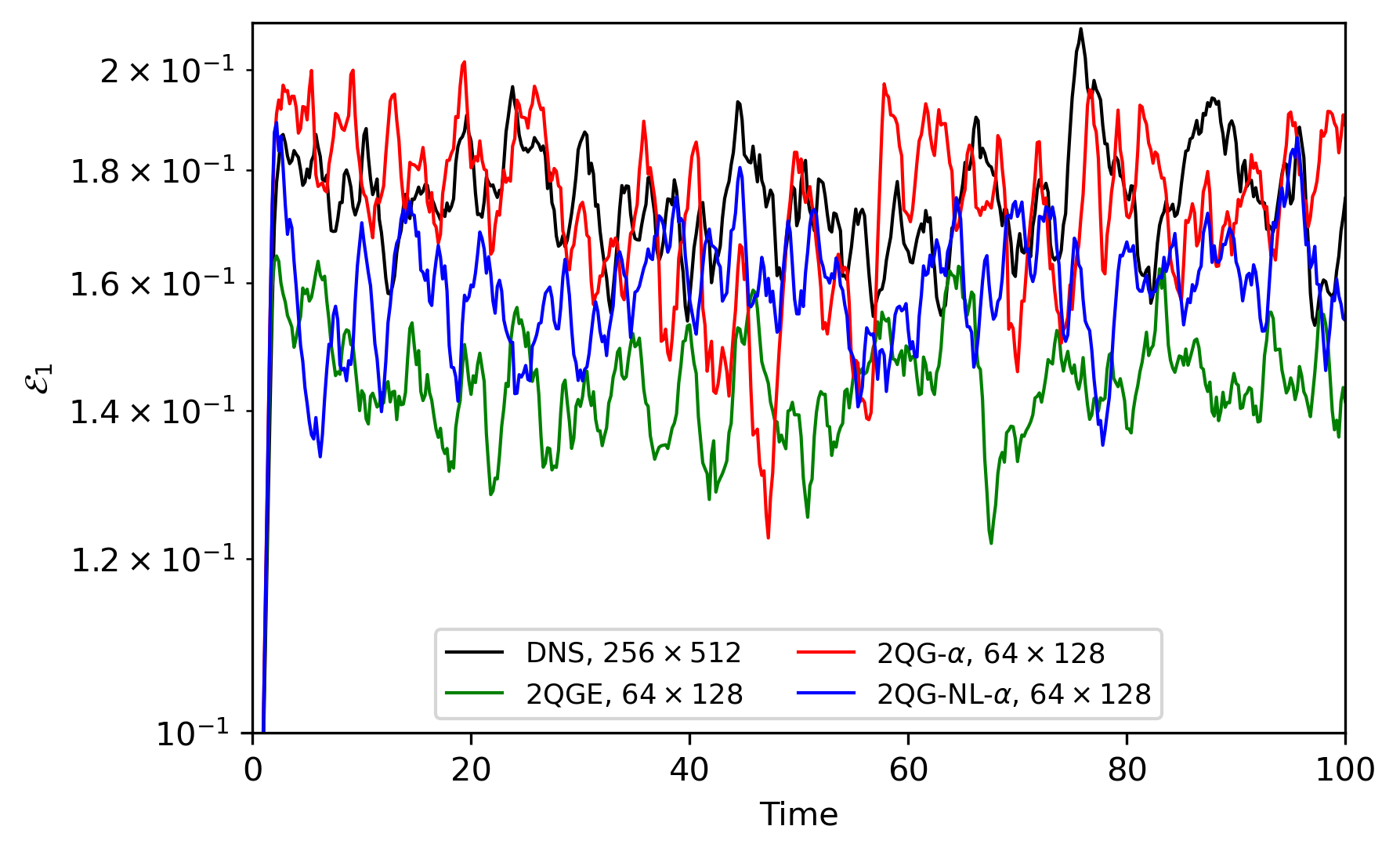}
     \end{subfigure}
    \begin{subfigure}{0.48\textwidth}
         \centering
         \includegraphics[width=\textwidth]{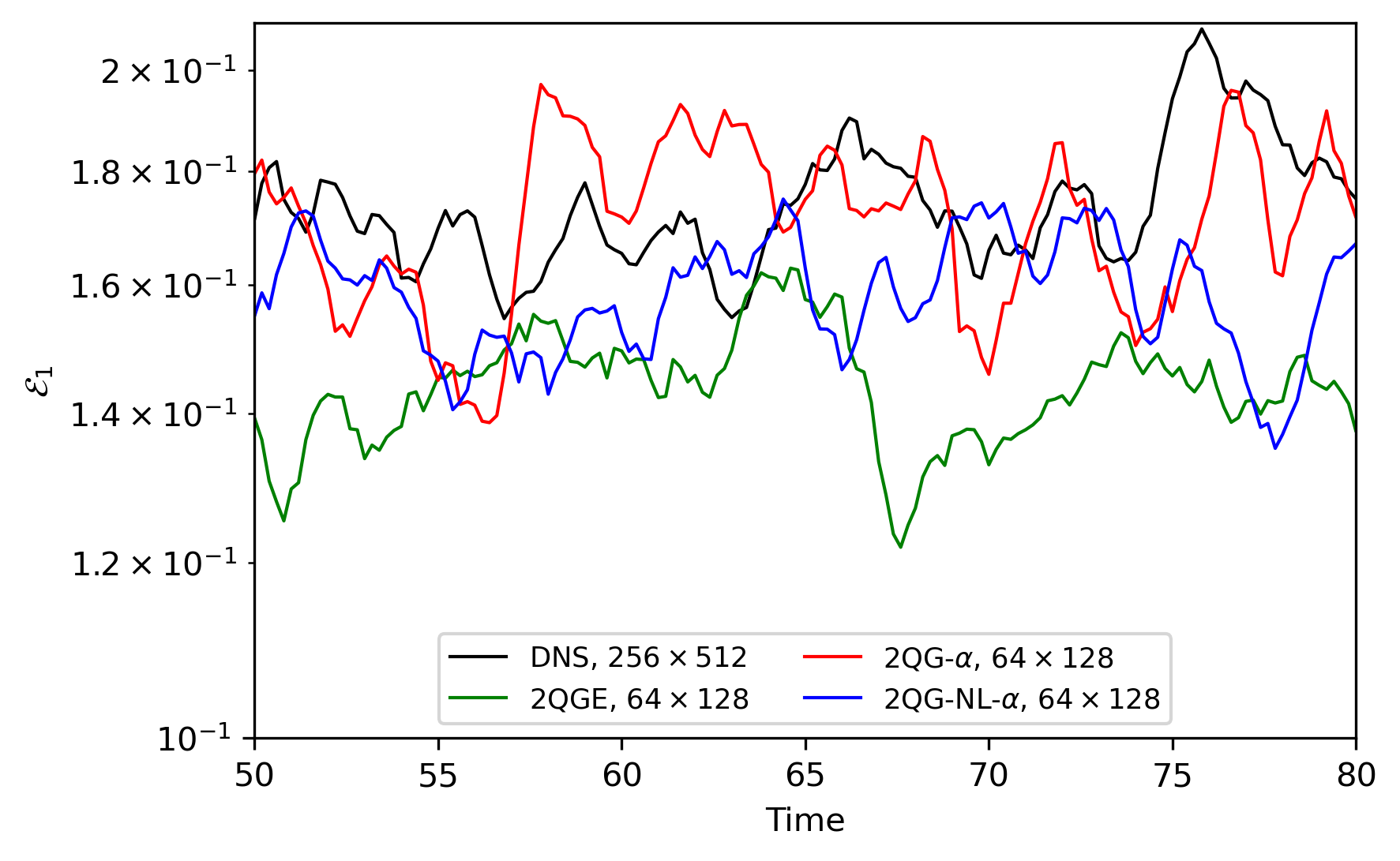}
     \end{subfigure}
    \begin{subfigure}{0.48\textwidth}
         \centering
         \includegraphics[width=\textwidth]{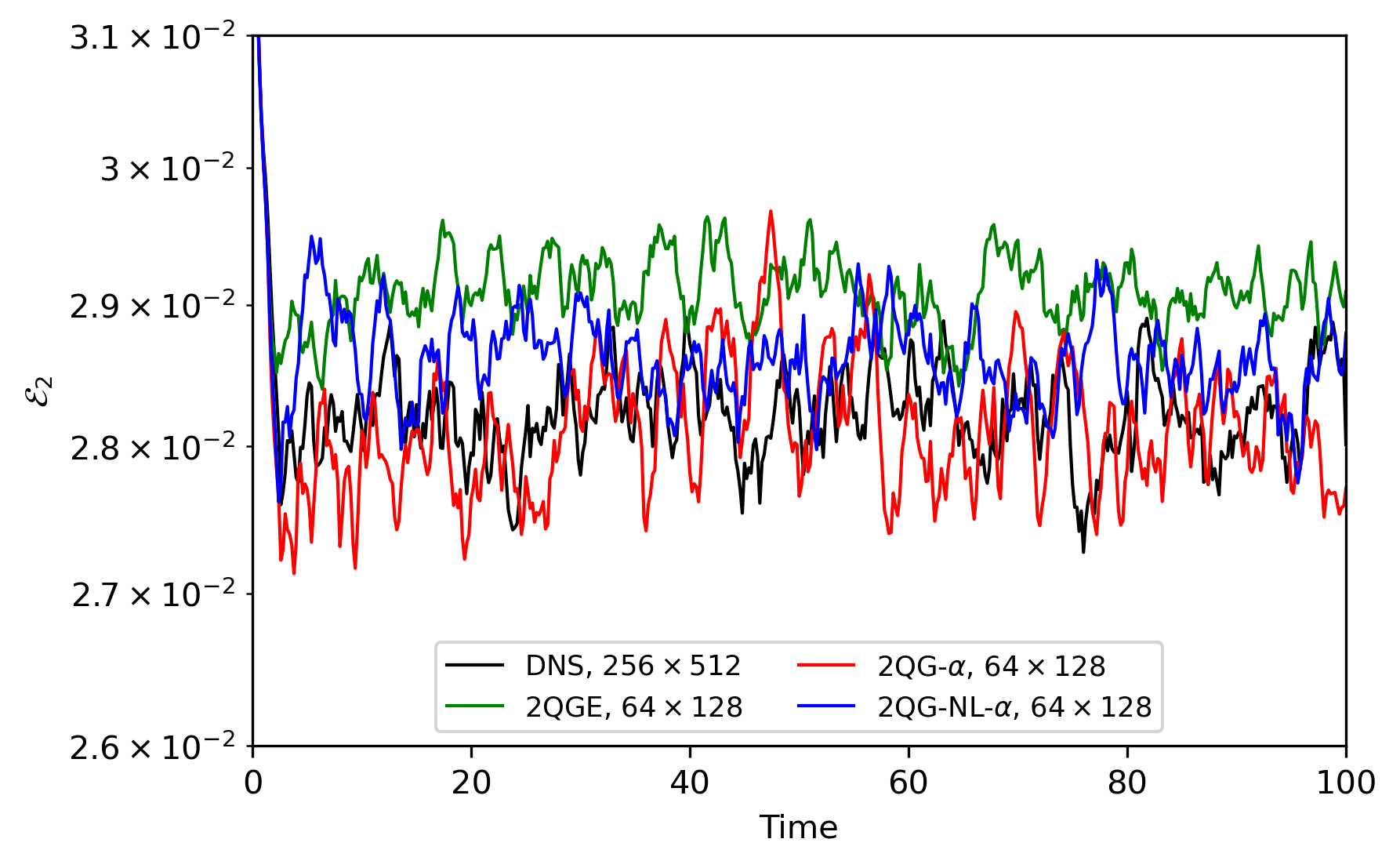}
     \end{subfigure}
    \begin{subfigure}{0.48\textwidth}
         \centering
         \includegraphics[width=\textwidth]{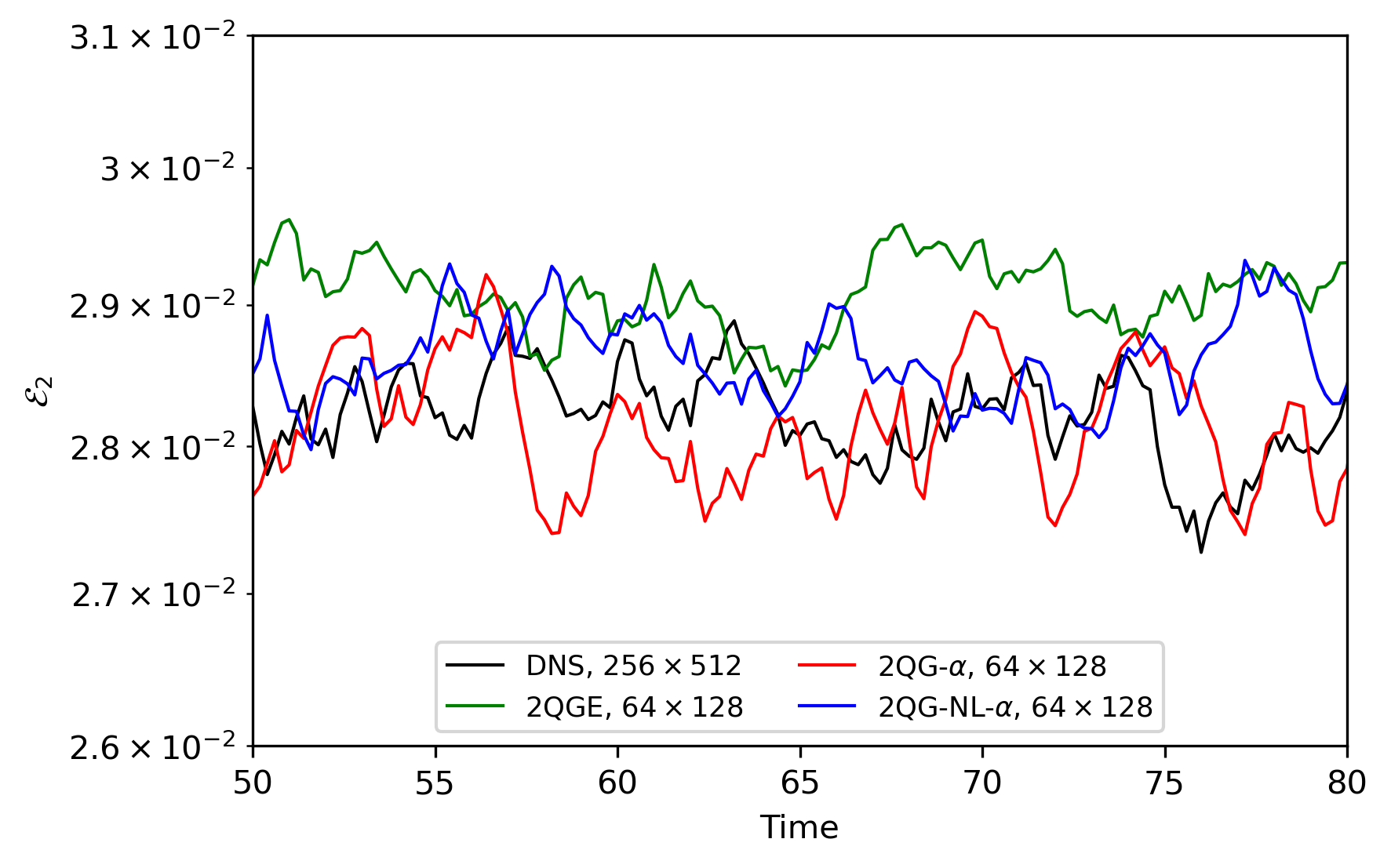}
     \end{subfigure}
\caption{Case 2: Time evolution of the enstrophy of the top layer $\mathcal{E}_1$ (top left) 
and the enstrophy of the bottom layer $\mathcal{E}_2$ (bottom left) 
for the DNS (mesh $256\times 512$), and 
2QGE model, 2QG-$\alpha$ model and 
2QG-NL-$\alpha$ model with mesh $64 \times 128$. 
The panels on the right are zoomed-in views
of the panels on the left.}
\label{fig:case2-ens}
\end{figure}

Fig.~\ref{fig:case2-ens} shows the computed enstrophy \eqref{eq:entrophy} of the two layers 
for the DNS, and the 2QGE, 2QG-$\alpha$, 2QG-NL-$\alpha$ models with a coarse mesh 
($64\times128$).
Similarly to case 1, the 2QGE with the coarse mesh
underestimates the enstrophy of the top layer 
and overestimates the enstrophy of the bottom layer.
For the given mesh, both the 2QG-$\alpha$ and 2QG-NL-$\alpha$ models improve the reconstruction of the time evolution of $\mathcal{E}_1$ and $\mathcal{E}_2$. However, it is hard to judge which does better. 
For this reason, we report in Table \ref{tab:exp2-ens-64x128}
the extreme values and average values of the enstrophies of both layers, together with the $L^2$ errors with respect to the DNS for the simulations with mesh $64\times128$. We note that
the 2QG-$\alpha$ and 2QG-NL-$\alpha$ models give comparable
values, with the 2QG-$\alpha$ slightly more accurate.

\begin{table}[h!]
    \centering
    \begin{tabular}{|c|c|c|c|c|c|c|}
       \hline
       \multirow{2}{*}{}  & \multirow{2}{*}{} & \multirow{2}{*}{DNS} & \multicolumn{3}{c|}{$64\times128$} \\ \cline{4-6} 
       \multirow{2}{*}{} & \multirow{2}{*}{} &  & 2QGE & 2QG-$\alpha$ & 2QG-NL-$\alpha$ \\ \hline
       \multirow{3}{*}{$\mathcal{E}_1$} & average & 1.732E-01 & 1.433E-01 & {1.718E-01} & 1.583E-01 \\ \cline{2-6}
       & max & 2.087E-01 & 1.646E-01 & {2.016E-01} & 1.892E-01 \\ \cline{2-6}
       & L2 error & - & 7.275E-01 & {3.750E-01} & 4.610E-01 \\ \hline
       \multirow{3}{*}{$\mathcal{E}_2$} & average & 2.826E-02  & 2.912E-02 & {2.816E-02} & 2.865E-02 \\ \cline{2-6}
       & min & 2.728E-02 & 2.840E-02 & {2.713E-02} & 2.762E-02 \\ \cline{2-6}
       & L2 error & - & 2.130E-02 & {1.180E-02} & 1.280E-02 \\ \hline
    \end{tabular}
    \caption{Case 2: extrema of $\mathcal{E}_1$ and $\mathcal{E}_2$ for the
    DNS (mesh $256\times512$) and the 2QGE, 2QG-$\alpha$, 2QG-NL-$\alpha$ models with coarse mesh
    $64 \times 128$.
    For the 2QGE, 2QG-$\alpha$, 2QG-NL-$\alpha$ models, we also
    report the associated $L^2$ errors.
    }
    \label{tab:exp2-ens-64x128}
\end{table}

We conclude with a comment on the computational time
required by the filtering techniques. 
Table \ref{tab:cpu_time_2} reports the CPU times 
for the DNS and the 2QGE, 2QG-$\alpha$, and 2QG-NL-$\alpha$
with coarser meshes $64\times128$ and $32\times64$, 
using the same machine as in case 1.
As one would expect, the computational times for the DNS and the simulations using mesh $32\times64$ are comparable to  those of case 1, obtaining a speed-up of 100-120 over the DNS when using filtering with mesh $32\times64$.
When using mesh $64\times128$, the computational times 
is 6 hours and 59 minutes for the 2QG-$\alpha$ model (30 times faster than the DNS) and 7 hours and 21 minutes for the 2QG-NL-$\alpha$ (28 times faster than the DNS). Thus, 
even in the case of more complex flow patterns, 
the filtering techniques, which allow to reconstruct
the time-averaged flow in an accurate manner, 
are considerably cheaper than the DNS.
Again, we see that the price to pay for the additional accuracy
of the nonlinear filter is small when compared to the
price of the linear filter. 

\begin{table}[htb!]
    \centering
    \begin{tabular}{lllc}
        \hline 
        Model & Mesh & CPU Time & Speed Up Factor \\
        \hline
        DNS & $256\times512$ & 8d 16h 27m & - \\ 
        2QGE & $64\times128$ & 5h 52m & - \\ 
        2QG-$\alpha$ & $64\times128$ & 6h 59m & 30 \\ 
        2QG-NL-$\alpha$ & $64\times128$ & 7h 21m & 28 \\ 
        2QGE & $32\times64$ & 1h 20m & - \\ 
        2QG-$\alpha$ & $32\times64$ & 1h 39m & 126 \\ 
        2QG-NL-$\alpha$ & $32\times64$ &  1h 54m & 109 \\ 
        \hline
    \end{tabular}
    \caption{
    Case 2: computational time required by the DNS and
    the 2QGE, 2QG-$\alpha$, and 2QG-NL-$\alpha$ models with coarser meshes $64\times128$ and $32\times64$. For the algorithms with filtering, the
    speed up factor with respect to the DNS is also reported.}
    \label{tab:cpu_time_2}
\end{table}

\section{Conclusions and future perspective}\label{sec:concl}

We applied linear and nonlinear filtering stabilization to the 2QGE with the goal of enabling significant computational time savings, through the use of coarser meshes, without sacrificing accuracy. The linear filtering technique, called 2QG-$\alpha$, introduces the same amount of regularization everywhere in the domain. The nonlinear filtering method, called 2QG-NL-$\alpha$, improves the methodology with the introduction of an indicator function that {identifies} where and how much regularization is needed. 
For the efficient implementation of the filtering techniques, we propose a segregated algorithm and space discretization by a Finite Volume method.

After validating the 2QGE solver and testing its robustness, we performed a computational study using a double-gyre wind forcing benchmark. The solutions given by the 2QG-$\alpha$ and 2QG-NL-$\alpha$ models
with coarse meshes
are compared to the solution given by a DNS, i.e., a simulation that resolves the Munk scale.
Our study indicates that, for a given coarse mesh, the 2QG-NL-$\alpha$ provides a more accurate reconstruction of large scale structures of the ocean dynamics compared to 2QG-$\alpha$. Both the 2QG-$\alpha$ and 2QG-NL-$\alpha$
models also outperform the 2QGE model 
in capturing the evolution of enstrophy using coarse meshes.

The indicator function used in this paper identifies the regions
of the domain with strong gradients of the potential vorticities.
This could be seen as a Leith-like method for the
2QGE. A possible extension of the work presented in this paper 
is the identification of more sophisticated indicator functions, 
for example by combining a filter with an approximate deconvolution
operator \cite{BQV,Girfoglio2019}. Another possible extension is to filter the stream-functions like in \cite{Holm2005}, 
to see if that outperforms our approach, which filters the potential vorticites.

\section*{Aknowledgements}
We acknowledge the support provided by 
PRIN “FaReX - Full and Reduced order modelling of coupled systems: focus on non-matching methods and automatic learning” project, PNRR NGE iNEST “Interconnected Nord-Est Innovation Ecosystem” project, INdAM-GNCS 2019–2020 projects and PON “Research and Innovation on Green related issues” FSE REACT-EU 2021 project. This work was also partially supported by the U.S. National Science Foundation through Grant No. DMS-1953535
(PI A. Quaini).

\bibliography{QGE} 

\begin{thebibliography}{10}
\expandafter\ifx\csname url\endcsname\relax
  \def\url#1{\texttt{#1}}\fi
\expandafter\ifx\csname urlprefix\endcsname\relax\def\urlprefix{URL }\fi
\expandafter\ifx\csname href\endcsname\relax
  \def\href#1#2{#2} \def\path#1{#1}\fi

\bibitem{Majda_Wang_2006}
A.~Majda, X.~Wang, {Nonlinear Dynamics and Statistical Theories for Basic
  Geophysical Flows}, Cambridge University Press, 2006.

\bibitem{Vallis_2017}
G.~K. Vallis, {Atmospheric and Oceanic Fluid Dynamics: Fundamentals and
  Large-Scale Circulation}, 2nd Edition, Cambridge University Press, 2017.

\bibitem{Marshall1997}
J.~Marshall, C.~Hill, L.~Perelman, A.~Adcroft, Hydrostatic, quasi-hydrostatic,
  and nonhydrostatic ocean modeling, Journal of Geophysical Research: Oceans
  102~(C3) (1997) 5733--5752.
\newblock \href {http://dx.doi.org/https://doi.org/10.1029/96JC02776}
  {\path{doi:https://doi.org/10.1029/96JC02776}}.

\bibitem{Chassignet1998}
T.~M. Ozgokmen, E.~P. Chassignet, Emergence of inertial gyres in a two-layer
  quasigeostrophic ocean model, Journal of Physical Oceanography 28~(3) (1998)
  461 -- 484.
\newblock \href
  {http://dx.doi.org/10.1175/1520-0485(1998)028<0461:EOIGIA>2.0.CO;2}
  {\path{doi:10.1175/1520-0485(1998)028<0461:EOIGIA>2.0.CO;2}}.

\bibitem{Berloff1999}
P.~S. Berloff, J.~C. McWilliams, Large-scale, low-frequency variability in
  wind-driven ocean gyres, Journal of Physical Oceanography 29~(8) (1999) 1925
  -- 1949.
\newblock \href
  {http://dx.doi.org/10.1175/1520-0485(1999)029<1925:LSLFVI>2.0.CO;2}
  {\path{doi:10.1175/1520-0485(1999)029<1925:LSLFVI>2.0.CO;2}}.

\bibitem{DiBattista2001}
M.~DiBattista, A.~Majda, Equilibrium statistical predictions for baroclinic
  vortices: The role of angular momentum, Theoret. Comput. Fluid Dynamics 14
  (2001) 293–322.

\bibitem{BERLOFF_KAMENKOVICH_PEDLOSKY_2009}
P.~Berloff, I.~Kamenovich, J.~Pedlosky, A mechanism of formation of multiple
  zonal jets in the oceans, Journal of Fluid Mechanics 628 (2009) 395–425.
\newblock \href {http://dx.doi.org/10.1017/S0022112009006375}
  {\path{doi:10.1017/S0022112009006375}}.

\bibitem{Bryan1963}
K.~Bryan, A numerical investigation of a nonlinear model of a wind-driven
  ocean, Journal of Atmospheric Sciences 20~(6) (1963) 594 -- 606.
\newblock \href
  {http://dx.doi.org/10.1175/1520-0469(1963)020<0594:ANIOAN>2.0.CO;2}
  {\path{doi:10.1175/1520-0469(1963)020<0594:ANIOAN>2.0.CO;2}}.

\bibitem{Gates1968}
W.~L. Gates, {A Numerical Study of Transient Rossby Waves in a Wind-Driven
  Homogeneous Ocean}, Journal of Atmospheric Sciences 25~(1) (1968) 3 -- 22.
\newblock \href
  {http://dx.doi.org/10.1175/1520-0469(1968)025<0003:ANSOTR>2.0.CO;2}
  {\path{doi:10.1175/1520-0469(1968)025<0003:ANSOTR>2.0.CO;2}}.

\bibitem{Holland1975}
W.~R. Holland, L.~B. Lin, On the generation of mesoscale eddies and their
  contribution to the oceanicgeneral circulation. ii. a parameter study,
  Journal of Physical Oceanography 5~(4) (1975) 658 -- 669.
\newblock \href
  {http://dx.doi.org/10.1175/1520-0485(1975)005<0658:OTGOME>2.0.CO;2}
  {\path{doi:10.1175/1520-0485(1975)005<0658:OTGOME>2.0.CO;2}}.

\bibitem{Tanaka2010}
Y.~Tanaka, K.~Akitomo, Alternating zonal flows in a two-layer wind-driven
  ocean, J. Oceanogr. 66 (2010) 475–487.

\bibitem{Holm2005}
D.~Holm, B.~Wingate, Baroclinic instabilities of the two-layer quasigeostrophic
  alpha model, Journal of Physical Oceanography 35 (2005) 1287--1296.

\bibitem{San2012}
O.~San, A.~Staples, T.~Iliescu, Approximate deconvolution large eddy simulation
  of a stratified two-layer quasigeostrophic ocean model, Ocean Modelling 63
  (2012) 1--20.

\bibitem{San2011}
O.~San, A.~Staples, Z.~Wang, T.~Iliescu, Approximate deconvolution large eddy
  simulation of a barotropic ocean circulation model, Ocean Modelling 40 (2011)
  120--132.

\bibitem{San_IJMCE2013}
O.~San, A.~E. Staples, An efficient coarse grid projection method for
  quasigeostrophic models of large-scale ocean circulation, International
  Journal for Multiscale Computational Engineering 11 (2013) 463--495.

\bibitem{Maulik2016}
R.~Maulik, O.~San, Dynamic modeling of the horizontal eddy viscosity
  coefficient for quasigeostrophic ocean circulation problems, Journal of Ocean
  Engineering and Science 1 (2016) 300--324.

\bibitem{Berloff2021}
P.~Berloff, E.~Ryzhov, I.~Shevchenko, On dynamically unresolved oceanic
  mesoscale motions, Journal of Fluid Mechanics 920 (2021) A41.

\bibitem{Salmon1978}
R.~Salmon, Two-layer quasi-geostrophic turbulence in a simple special case,
  Geophysical \& Astrophysical Fluid Dynamics 10 (1978) 25--52.

\bibitem{Medjo2000}
T.~T. Medjo, Numerical simulations of a two-layer quasi-geostrophic equation of
  the ocean, SIAM Journal of Numerical Analysis 37 (2000) 2005--2022.

\bibitem{Fandry1984}
C.~B. Fandry, L.~M. Leslie, A two-layer quasi-geostrophic model of summer
  trough formation in the australian subtropical easterlies, Journal of the
  Atmospheric Sciences 41 (1984) 807--818.

\bibitem{Mu1994}
M.~Mu, Z.~Qingcun, T.~G. Shepherd, L.~Yongming, Nonlinear stability of
  multilayer quasi-geostrophic flow, Journal of Fluid Mechanics 264 (1994)
  165--184.

\bibitem{Girfoglio_JCAM2023}
M.~Girfoglio, A.~Quaini, G.~Rozza, A novel large eddy simulation model for the
  quasi-geostrophic equations in a finite volume setting, Journal of
  Computational and Applied Mathematics 418 (2023) 114656.

\bibitem{Monteiro2015}
I.~Monteiro, C.~Manica, L.~Rebholz, Numerical study of a regularized barotropic
  vorticity model of geophysical flow, Numerical Methods for Partial
  Differential Equations 31 (2015) 1492--1514.

\bibitem{Nadiga2001}
B.~Nadiga, L.~Margolin, Dispersive-dissipative eddy parameterization in a
  barotropic model, Journal of Physical Oceanography 31 (2001) 2525--2531.

\bibitem{Holm2003}
D.~Holm, B.~Nadiga, Modeling mesoscale turbulence in the barotropic double-gyre
  circulation, Journal of Physical Oceanography 33 (2003) 2355--2365.

\bibitem{Monteiro2014}
I.~Monteiro, C.~Carolina, Improving numerical accuracy in a regularized
  barotropic vorticity model of geophysical flow, International Journal of
  Numerical Analysis and Modeling, Series B 5 (2014) 317--338.

\bibitem{Girfoglio2019}
M.~Girfoglio, A.~Quaini, G.~Rozza, A finite volume approximation of the
  {Navier-Stokes} equations with nonlinear filtering stabilization, Computers
  \& Fluids 187 (2019) 27--45.

\bibitem{GirfoglioPSIZETA}
M.~Girfoglio, A.~Quaini, G.~Rozza, {A POD-Galerkin reduced order model for the
  Navier-Stokes equations in stream function-vorticity formulation}, Computers
  \& Fluids 244 (2022) 105536.

\bibitem{GEA}
{GEA - Geophysical and Environmental Applications},
  \url{https://github.com/GEA-Geophysical-and-Environmental-Apps/GEA}.

\bibitem{GirfoglioFVCA10}
M.~Girfoglio, A.~Quaini, G.~Rozza, {GEA: A New Finite Volume-Based Open Source
  Code for the Numerical Simulation of Atmospheric and Ocean Flows}, in: Finite
  Volumes for Complex Applications X—Volume 2, Hyperbolic and Related
  Problems, 2023, pp. 151--159.

\bibitem{Weller1998}
H.~G. Weller, G.~Tabor, H.~Jasak, C.~Fureby, A tensorial approach to
  computational continuum mechanics using object-oriented techniques, Computers
  in physics 12 (1998) 620--631.

\bibitem{Girfoglio_AIP2023}
M.~Girfoglio, A.~Quaini, G.~Rozza, Validation of an
  {OpenFOAM\textsuperscript{\textregistered{}}-based} solver for the {Euler
  equations} with benchmarks for mesoscale atmospheric modeling, AIP Advances
  13 (2023) 055024.

\bibitem{Clinco_2023}
N.~Clinco, M.~Girfoglio, A.~Quaini, G.~Rozza, Filter stabilization for the
  mildly compressible euler equations with application to atmosphere dynamics
  simulations, Computers \& Fluids 266 (2023) 106057.

\bibitem{Hajisharifi2024}
A.~Hajisharifi, M.~Girfoglio, A.~Quaini, G.~Rozza, A comparison of data-driven
  reduced order models for the simulation of mesoscale atmospheric flow, Finite
  Elements in Analysis and Design 228 (2024) 104050.

\bibitem{QGE-review}
C.~Mou, Z.~Wang, D.~Wells, X.~Xie, T.~Iliescu, Reduced order models for the
  quasi-geostrophic equations: A brief survey, Fluids 6 (2021) 16.
\newblock \href {http://dx.doi.org/10.3390/fluids6010016}
  {\path{doi:10.3390/fluids6010016}}.

\bibitem{Greathbatch2000}
R.~Greatbatch, B.~Nadiga, Four-gyre circulation in a barotropic model with
  double-gyre wind forcing, Journal of Physical Oceanography 30 (2000)
  1461--1471.

\bibitem{San2015}
O.~San, T.~Iliescu, A stabilized proper orthogonal decomposition reduced-order
  model for large scale quasigeostrophic ocean circulation, Advances in
  Computational Mathematics 41 (2015) 1289--1319.

\bibitem{Girfoglio2023}
M.~Girfoglio, A.~Quaini, G.~Rozza, A linear filter regularization for
  {POD-based} reduced-order models of the quasi-geostrophic equations, Comptes
  Rendus. M\'ecanique 351 (2023) 1--21.

\bibitem{BQV}
L.~Bertagna, A.~Quaini, A.~Veneziani, {Deconvolution-based nonlinear filtering
  for incompressible flows at moderately large Reynolds numbers}, Int. J.
  Numer. Meth. Fluids 81~(8) (2016) 463--488.

\end{thebibliography}

\end{document}